\documentclass[preprint,12pt]{elsarticle}

\usepackage{graphicx}
\usepackage{amssymb}
\usepackage{amsmath}
\usepackage{amsthm}

\usepackage{subfig}

\usepackage[pagewise]{lineno}

\newtheorem{theorem}{Theorem}
\newtheorem{remark}[theorem]{Remark}
\newtheorem{conjecture}[theorem]{Conjecture}
\newtheorem*{acknowledgement}{Acknowledgement}

\journal{Computer Methods in Applied Mechanics and Engineering}

\begin{document} 
\begin{frontmatter}



\title{On the $p$-version of FEM in one dimension: The known and unknown features}


\author[ib]{Ivo Babu\v{s}ka}
\address[ib]{
The University of Texas at Austin\\ 
Department of Aerospace Engineering and 
   Engineering Mechanics\\ 
105 W. Dean Keeton Street, SHC 328\\ 
Austin, TX 78712, USA
}
\ead{babuska@ices.utexas.edu}

\author[hh]{Harri Hakula\fnref{support}}
\fntext[support]{H. Hakula: Supported by T J Oden Fellowship, 
ICES University of Texas at Austin.}

\address[hh]{Aalto University\\
Department of Mathematics and System Analysis\\
P.O. Box 11100\\
FI--00076 Aalto, Finland
}
\ead{Harri.Hakula@aalto.fi}

\begin{abstract}
The paper analyses the convergence of the $p$-version of the FEM when the solution is 
piecewise analytic function. It focuses on pointwise convergence of the gradient. 
It shows that at boundary the rate is different than inside the element and there is 
a Gibbs phenomenon in the neighborhood of the point where the solution is not analytic. 
The major result is a conjecture written in the form of a theorem. 
The conjecture is based on careful numerical computations. 
The known theoretical results are stated.  
\end{abstract}


\end{frontmatter}


\section{Introduction}

There are two basic versions of the Finite Element Method; The $h$-version
and the $p$-version. In addition the versions can be combined into
the $hp$-version. 
The first papers on the $p$-version and the $hp$-version were \cite{c1,c2}. 
In these papers the major basic theorems on the errors 
in the energy norm were established.
The $h$-version with its roots in engineering
problems is much older, however. The
\textit{mathematical} theory was developed in the early 1970's. We refer to  \cite{c3,c4} 
where the history up to 1970 is presented. For the basic (older)
monographs on the theory of the $h$-version we refer for example to  \cite{c5,c6}.
Since then the theory of the $h$-version has been developed in many 
different ways
and there is a large body of literature on this subject. 
Especially we would like to
mention the theory of the error estimates in norms other than the energy
norm. For the local behavior of the $h$-version and pointwise estimates we
refer to the excellent survey \cite{c7} and literature mentioned therein. All
theoretical results for the $h$-version are independent of the degree of the
elements.

Situation is different for the $p$-version. The domain is divided into
elements as in the $h$-version but the number of elements is fixed and the
convergence is obtained by increasing the degrees $p$ either 
uniformly or adaptively.
In the $hp$-version the convergence is obtained by the simultaneous 
refinement of
the finite element mesh. 
For computation with the $hp$-version we refer
to the books \cite{bc1,bc2,c20,c21} and citations therein, and the commercial software
STRESSCHECK. 
The existing theory addresses only the error
estimates in the energy and $L^{2}$-type norms. 
Thus, general theory of the local behavior is not available. 

In one dimension the $p$-version is closely related to the expansion of
functions in the Legendre series. Indeed, the subject of approximation and interpolation
of a function
by polynomials is a classical one. We refer here to the excellent books of
the approximation theory and orthogonal polynomials in one
dimension \cite{c9,c10}.
There is large literature on orthogonal polynomials, see e.g, \cite{c11},
and polynomials in general, e.g, \cite{c12}. 
For the history of the numerical problems we refer
to the second chapter of \cite{c13} and to \cite{c14}. For more recent books on the 
$p$-version also called the spectral method we refer to \cite{c8}, \cite{c15}, 
and the recent book \cite{cm2} with citations.

Nevertheless, despite all available results there are still open theoretical
problems about the local behavior and the pointwise convergence of the 
$p$-version in one and higher dimensions. 
This paper addresses the behavior of the $p$-version when the solution is piecewise analytic. 
Such class of
solutions is typical in the applications. Here the focus is on one dimensional
problems. The paper is a survey of known results and it formulates new
theorems as conjectures based on the detailed computations. 
It also serves as an introduction to the analysis of the $p$-version in two and three dimensions 
which will be addressed in the next papers. There are various similarities between the behavior 
of the $p$-version in one and higher dimensions as well significant differences. 
For example in one dimension there was no pollution from one element into the neighboring ones. 
This is important when the elements are refined in the neighborhood of a singularity. 
For some results in this direction we refer to \cite[p192]{bc1}  and \cite[p190]{bc2}.
In higher dimensions solving the system of linear equations became nontrivial task. 
Recently various iterative methods were analyzed. 
Recently in \cite{bc3} a multigrid method for solving the problem in one three dimensional 
cube with $p =16$. 
Nevertheless this paper and the next ones are concentrating only on the problem of 
the convergence, Gibbs phenomenon, and element boundary behavior among others.

Through a series of numerical experiments with highly accurate results we arrive at
the main theoretical result of this paper, the theorem (conjecture) on the
error of the Legendre expansion of piecewise analytic functions,
\[
  u(x)=\sum_{i=1}^{m}c_{i}\mid x-a_{i}\mid ^{\beta _{i}}+\ v(x),\mid
a_{i}\mid <1,\beta _{i}>-1
\]
and $v(x)$ {analytic function} on [-1,1].
This conjecture is fully stated in Section~\ref{sec:hypo} below. The stated convergence rates 
are optimal in the sense that they cannot be improved.
We emphasize that the existing mathematical theory does not cover all
parts of the conjecture and thus hopefully will lead to more
refined mathematical analysis in the future.

The paper is organized as follows. In the Section 2 we define the one dimensional
problem and its solution $u$. Section 3 shows that the $p$-version in this
simple setting is the Legendre expansion of $u^{\prime}$. In the
Section 4 we analyze the Legendre expansion of $u^{\prime}$, the error of its
partial sum, and comment on the known and not available results. Section
5 analyzes the Legendre expansion of the solution $u$. 
Section 6 addresses the $p$-version of the
solution which is the partial Legendre sum with a constraint at the boundary
points. Section 7 generalizes the problem addressed in the Section 5.
In Section 8 we summarize the result and formulate the theorem
(conjectures) on the behavior of the $p$-version when the
solution is piecewise analytic function.
\section{The Problem}
Consider the boundary value problem
\begin{equation}\label{eq:2.1}
-u^{\prime \prime }(x)=\delta_{a}(x),\ x \in I=(-1,1),
\end{equation}
with the boundary conditions
\begin{equation}\label{eq:2.1b}
u(\pm 1)=0.
\end{equation}
Here\ $\delta _{a}(x)$ is the Dirac function acting at $x=a.$

Obviously
\begin{equation}\label{eq:2.2}
u(x)= 
\left\{ \begin{array}{rl}
 c\,(x+1), &\mbox{ if $-1<x<a$}, \\
 c\,(a+1)+(1+c)(x-a), &\mbox{ if $a\leq x<1$},
\end{array} \right.
\end{equation}
where $c=-\frac{1-a}{2}$.
The solution $u\notin H_{0}^{3/2}(I)$, but $u\in $ $H_{0}^{1}(I)\cap B_{2,\infty }^{3/2}(I)$
where $B_{2,\infty }^{3/2}(I)$ is the Besov space.

The weak solution $u\in H_{0}^{1}(I)$ satisfies
\begin{equation}
B(u,v)=\int_{-1}^{1}u^{\prime }v^{\prime }dx=v(a),\forall v\in H_{0}^{1}(I).
\end{equation}
Let us denote by $\mathcal{E}_{0}(I)$\ $=H_{0}^{1}(I)$ the energy space with
the product $(u,v)_{\mathcal{E}_{{}}}.$

\section{The $h$- and $p$-versions of the Finite Element Method}
\label{sec:problem}

Let $\mathcal{M}$ be a uniform FEM mesh with the nodes 
$x_{i}^{(h)}=ih-1,i=0,1,\ldots,n, n=\frac{2}{h}$, and the (open) elements 
$\omega_{i}^{(h)}=(x_{i-1}^{(h)},x_{i}^{(h)})$. Let us assume that 
point $a\in \omega_{s}^{(h)}$.
Further let $S_{p}^{(h)}(I)\subset \mathcal{E}_{0}$ be the space of the
continuous piecewise continuous polynomials of degree $p$ on every element 
$\omega _{i}^{(h)}$ and $u_{p}^{(h)}\in S_{p}^{(h)}$, the $hp$-FEM solution of the 
Problem \ref{eq:2.1}, 
which satisfies
\begin{equation}
B(u_{p}^{(h)},v_{p})=v(a), \forall v_{p} \in S_{p}^{(h)}.
\end{equation}

From the classical FEM theory we have $u_{p}^{(h)}=u$ on every 
$\omega_{j}^{(h)},j\neq s$, i.e., on every element which does not contain the 
point $x=a$. Further denote $\widehat{u}^{(h)}\in S_{p}^{(h)}$ the piecewise
linear function which coincides with the exact solution $u$ at the nodes. 
Let us denote the difference $w_{{}}^{(h)}=u-\widehat{u}^{(h)}$. 
Then $w^{(h)}=0$ on every 
$\omega _{i}^{(h)},i\neq s$, and on the boundary of the element 
$\omega_{s}^{(h)}$containing the point $x=a$. 
Further let $w_{p}^{(h)}=u_{p}^{(h)}-\widehat{u}^{(h)}.$
Then $w_{p}^{(h)}=0$ on all $\omega _{j}^{(h)},j\neq s$,
and ($w_{p}^{(h)})^{\prime }$ is the $p$-partial sum of the Legendre
expansion of the function ($w^{(h)})^{\prime }$ on 
$\omega_{s}^{(h)}.$ Fixing the $p$ and letting $h\rightarrow 0$ we get 
the $h$-version of the FEM with the elements of order $p$. 
Fixing $h$ and letting $p\rightarrow \infty$ we get the $p$-version.

For small fixed $p$ the error $\varepsilon_{p}^{(h)}(x)=u(x)-u_{p}^{(h)}(x)$ 
can be easily estimated. For the $h$-version there are many available results. 
Not so for the $p$-version. Here
in this paper we shall be especially interested to study the pointwise
convergence for $\varepsilon_{p}^{(h)}$ and 
($\varepsilon_{p}^{(h)})^{\prime }$ for the $p$-version of the FEM. 

Thus, without any loss of the generality, we can
concentrate only on the element $\omega_{s}^{(h)}$ which contains $x=a$ 
and assume that $h=2$, i.e., to analyze the error of the Legendre
expansion when using  $p$-terms in the Legendre expansion.

\section{The Legendre expansion of the function $(u^{(h)})^{\prime}(x)$ with $h=2$}
\label{sec:legendre}

Here we write $\varepsilon_{p}$ resp. $u$ instead of $\varepsilon _{p}^{(h)}(x)$
resp. $u^{(h)}(x)$ and $x\in I=(-1,1)$.

Obviously we have
\begin{equation}\label{eq:4.1}
u^{\prime }(x)= 
\left\{ \begin{array}{rl}
c, &\mbox{ if $-1\leq x<a,$} \\
1/2+c, &\mbox{ if $x=a$,} \\
1+c, &\mbox{ if $a\leq x\leq 1$},
\end{array} \right.
\end{equation}
with $\int_{-1}^{1}u^{\prime }(x)dx=0$. From it we
get $c=\frac{a-1}{2}$.

Function $u^{\prime }(x)$ is discontinuous at $x=a$ and in (\ref{eq:4.1}) 
we define $u^{\prime }(a)=\frac{1}{2}(u^{\prime }(x+0)+u^{\prime }(x-0))$ 

Denoting the Legendre polynomials of degree $k$ with $P_k(x)$ we have
\begin{equation}
u^{\prime }(x)=\sum_{k=1}^{\infty }a_{k}P_{k}(x),
\end{equation}
where the coefficients
\begin{equation}\label{eq:4.3}
a_{k}=(k+\frac{1}{2})\int_{-1}^{1}u^{\prime
}(x)P_{k}(x)dx=\frac{1}{2}(P_{k-1}(a)-P_{k+1}(a)),k=1,2,\ldots
\end{equation}
For the error
\begin{equation}
\varepsilon _{p}^{\prime }(x)=u^{\prime }(x)- u_p^{\prime }(x)=u^{\prime }(x)-\sum_{k=1}^{p}a_{k}P_{k}(x).
\end{equation}
In the following we analyze $\varepsilon_{p}^{\prime}(x)$ as a function of 
$p$ for different $x.$

First let us mention the major theorem about the error of the Legendre
expansion.
\begin{theorem}[\cite{c16}]\label{thm:1}
Let $f(x)$ be a function of bounded variation on [-1,1].
Let $f(x)=\sum_{k=0}^{\infty }a_{k}(f)P_{k}(x)$ and 
$S_{p}(f,x)=\sum_{k=0}^{p}$\ $a_{k}(f)P_{k}(x)$ be 
the $p$-partial sum of the Legendre
series of $f$. One has $a_{k}(f)=(k+\frac{1}{2}) 
\int_{-1}^{1}f(t)P_{k}(t)dt.$ Then for $x\in (-1,1)$ and $p\geq 2$
	\begin{multline}\label{eq:4.4}
		\mid S_{p}(f,x)-\frac{1}{2}(f(x+0)-f(x-0))\mid
\leq \\ \frac{28}{p}(1-x^{2})^{-3/2}%
\sum_{k=1}^{p}V_{x-(1+x)/k}^{x+(1-x)/k}(g_{x})+ \\ \frac{1}{\pi p}%
(1-x^{2})^{-1}\mid f(x+0)-f(x-0)\mid,
		\end{multline}
		where
	\begin{equation}
		g_{x}(t)=
\left\{ \begin{array}{rl}
f(t)-f(x-0), &\mbox{ if $-1\leq t<x,$} \\
0, &\mbox{ if $t=x$} \\
f(t)-f(x+0), &\mbox{ if $x<t\leq 1$},
\end{array} \right.
		\end{equation}
and $V_{a}^{b}(g)$ is the total variation of $g$ on $[a,b]$.
\end{theorem}

\begin{remark}
  In the proof of the above theorem the inequality 
\begin{equation}
\mid P_{p}(x)\mid \leq \left(\frac{2}{\pi }%
\right)^{1/2}(1-x^{2}\ )^{-1/2}p^{-1/2}
\end{equation}
is used. Nevertheless later results \cite{c9} proved
\begin{equation}
\mid P_{p}(x)\mid \leq(1-x^{2})^{1/4}\sqrt{\frac{2}{\pi p}}.
\end{equation}
By this improved estimate we get the powers of $(1-x^{2})$ as $-5/4$ and $-1/2$
in the first and second terms, respectively.
\end{remark}

\begin{remark}
 With $a=0$ in (\ref{eq:4.1}) the following holds:
\begin{theorem}[\cite{cM}]\label{thm:2}
We have 
\[
  \mid \varepsilon _{p}^{\prime }(x)\mid \leq
  \frac{C}{p((1-x^{2})^{1/2}+p^{-1})^{1/2}}
\]for $2\delta <\mid x\mid \leq
1,\ 0<\delta <1/4,$ and $C$ independent of $p.$
 \end{theorem}
\end{remark}

The paper \cite{c17} is only a slight generalization of the Theorem \ref{thm:1} 
and the paper
\cite{c18} addresses Legendre expansion of functions which are analytic on [-1,1].
We are not aware of any other general theorems similar to
Theorem \ref{thm:1}.

Let us now apply the Theorem \ref{thm:1} to our example and consider the different cases
separately:
\begin{itemize}
\item[1.]$x=a.$ Then the second term in (\ref{eq:4.4}) is 
$(\pi p(1-a^{2}))^{-1}$ and the
first term is
\begin{equation}
\frac{28}{p}(1-a^{2})^{-3/2}%
\sum_{k=1}^{p}1+\frac{1}{\pi p}(1-x^{2})^{-1}=28(1-a^{2})^{-3/2}+\frac{1}{%
\pi p}(1-a^{2})^{-1}
\end{equation}

Theorem does not indicate the convergence, however,
we can compute the error directly. We have 
\begin{equation}\label{eq:4.9}
\mid \varepsilon_{p}^{\prime }(a)\mid =\mid u(a) - \frac{1}{2}%
 \sum_{k=1}^{p}(P_{k-1}(a)-P_{k+1}(a))P_{k}(a)=\mid \frac{1}{2}
P_{p+1}(a)P_{p}(a)\mid \leq \frac{1}{\pi p}(1-x^{2})^{-1/2}
\end{equation}

As we will see below $\parallel \varepsilon _{p}^{\prime }(x)\parallel
_{L^{\infty }(I)}\nrightarrow 0$ as $p\rightarrow \infty .$

\item[2.] $x\neq a,\pm 1.$ Then the second term in (\ref{eq:4.4}) is zero and we have
  \begin{equation}\label{eq:4.10}
\mid \varepsilon_{p}^{\prime}(x)\mid =\frac{28}{p}(1-x^{2})^{-3/2}\sum_{k=1}^{m}1=\frac{28}{p}(1-x^{2})^{-3/2}m,
\end{equation}
where for $x<a$, $m$ is the largest integer such that $x+(1-x)/k\geq 1+a$, and
for $x>a$, $m$ is the largest integer such that $x-(1+x)/k\geq (1-a)$.
\begin{remark}
  From (\ref{eq:4.10}) we see that the
error $\mid \varepsilon _{p}^{\prime }(x)\mid $ increases as $%
x\rightarrow a$ \textit{and for }$x=a-\frac{1}{p}$ we have $\mid
\varepsilon _{p}^{\prime }(x)\mid \sim 1$ which relates to the Gibbs
phenomenon discussed below. This is also the reason of exclusion of $x=0$ in the
Theorem~\ref{thm:2}.
\end{remark}
\item[3.]$x=\pm 1.$ Theorem 1 is not applicable to this case. 
Nevertheless in our
example we can compute the error directly. We get
\begin{equation}\label{eq:4.11}
\mid \varepsilon _{p}^{\prime }(\pm 1)\mid \leq \frac{1}{%
2}(\mid P_{p}(a)\mid +\mid P_{p+1}(a)\mid )\leq (1-a^{2})^{-1/4}\left(\frac{2}{%
\pi p}\right)^{1/2}
\end{equation}
\end{itemize}
Let us now show the computational results for the specific values $a=1/2$
and $c=-\frac{1}{4}$.
From the above theory we expect that
\begin{equation}\label{eq:4.12}
\mid \varepsilon _{p}^{\prime }(x)\mid \leq
C(x)\,p^{-\alpha (x)}
\end{equation}
with $\alpha $ and $C$ depending on $x$. We shall compute for $p\in (1,2200)$
and present the results in the loglog-scale. Notice that the upper limit 2200 is
chosen so that in double precision floating point arithmetic, the rounding errors
do not pollute the results in those cases where the series coefficients are known in
advance. In some instances, however, we are forced to compute using symbolic
computations in exact arithmetic. Naturally, this requires significant computing resources
in terms of time and memory.

\begin{figure}[ht]
	\centering
	\subfloat[Error at $x=-1$, $\alpha (-1)=1/2$ , $C(-1)=0.44194$.]
	{\label{fig:fig1a}\includegraphics[width=0.45\textwidth]{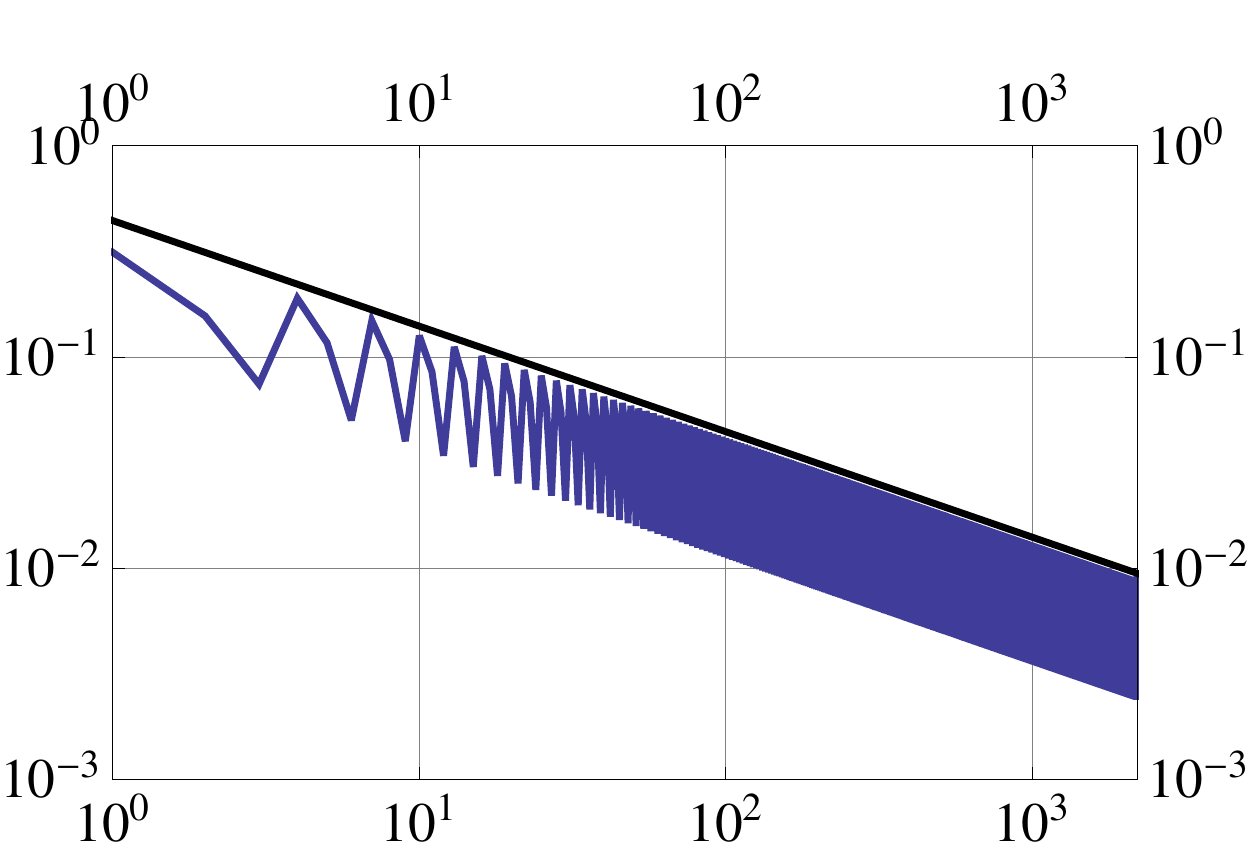}}\quad
	\subfloat[Error at $x=1$, $\alpha (1)=1/2$ , $C(1)=0.75$.]
	{\label{fig:fig1b}\includegraphics[width=0.45\textwidth]{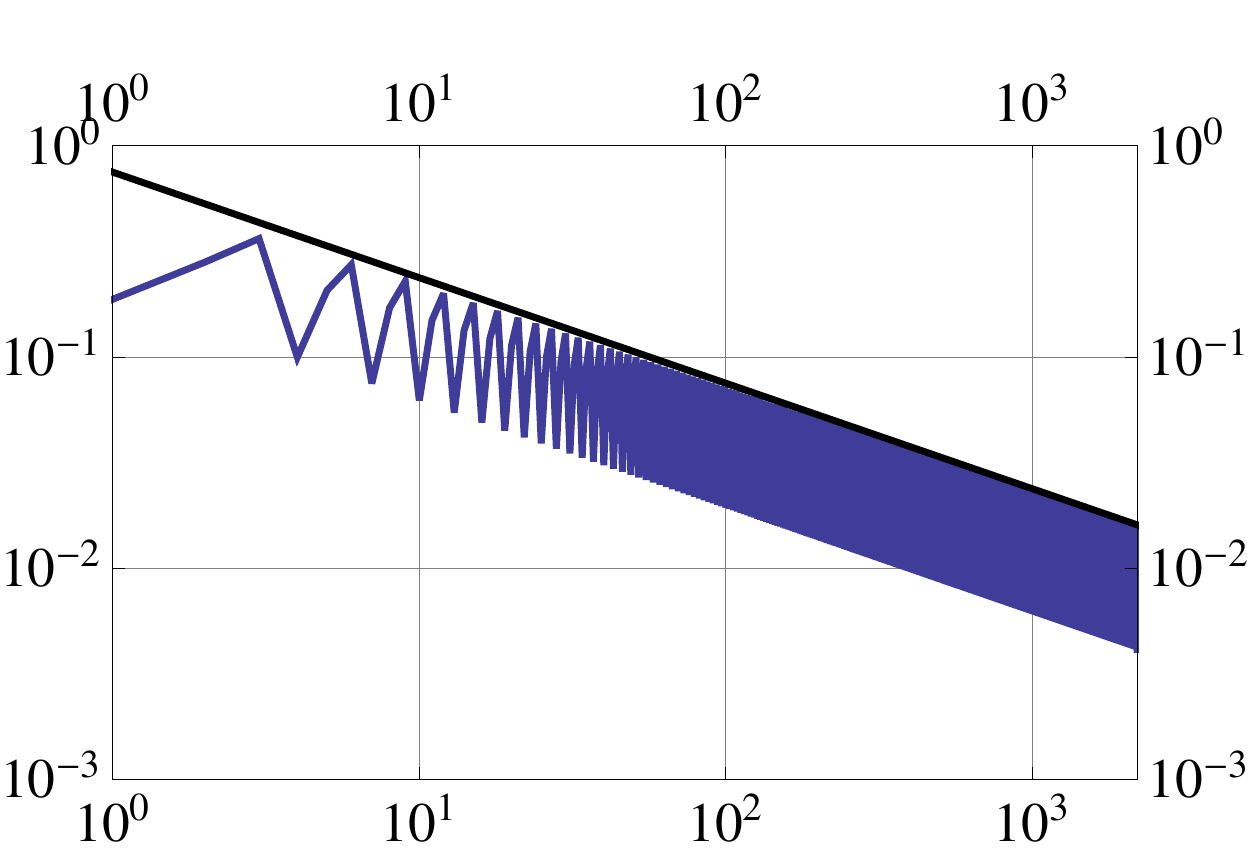}}
  \caption{$\varepsilon^{\prime }(x)$: Absolute value vs polynomial order, $p=1,\ldots,2200$.}
  \label{fig:fig1}
\end{figure}
Different behavior of the error for different values of $x$ will be seen.
Figures~\ref{fig:fig1a}--\ref{fig:fig1b} show the error for $x=\pm 1$. 
We see that the form (\ref{eq:4.12}) leads to
very accurate upper estimate with $\alpha (\pm 1)=1/2$, $C(-1)=0.44194$ for 
$x=-1$ and for $x=1$ we have $C(1)=0.75.$ From (\ref{eq:4.11}) we get $C=0.85738$.

There is slightly larger error for $x=1$ than for $x=-1$ because the
discontinuity is closer. The preasymptotic behavior is very similar as well.
It is interesting to notice that we also have $\mid \varepsilon ^{\prime }(\pm 1)\mid
\geq C^{\ast }p^{-1/2}$ with $C^{\ast }=0.15$.

\begin{figure}[ht]
	\centering
	{\includegraphics[width=0.45\textwidth]{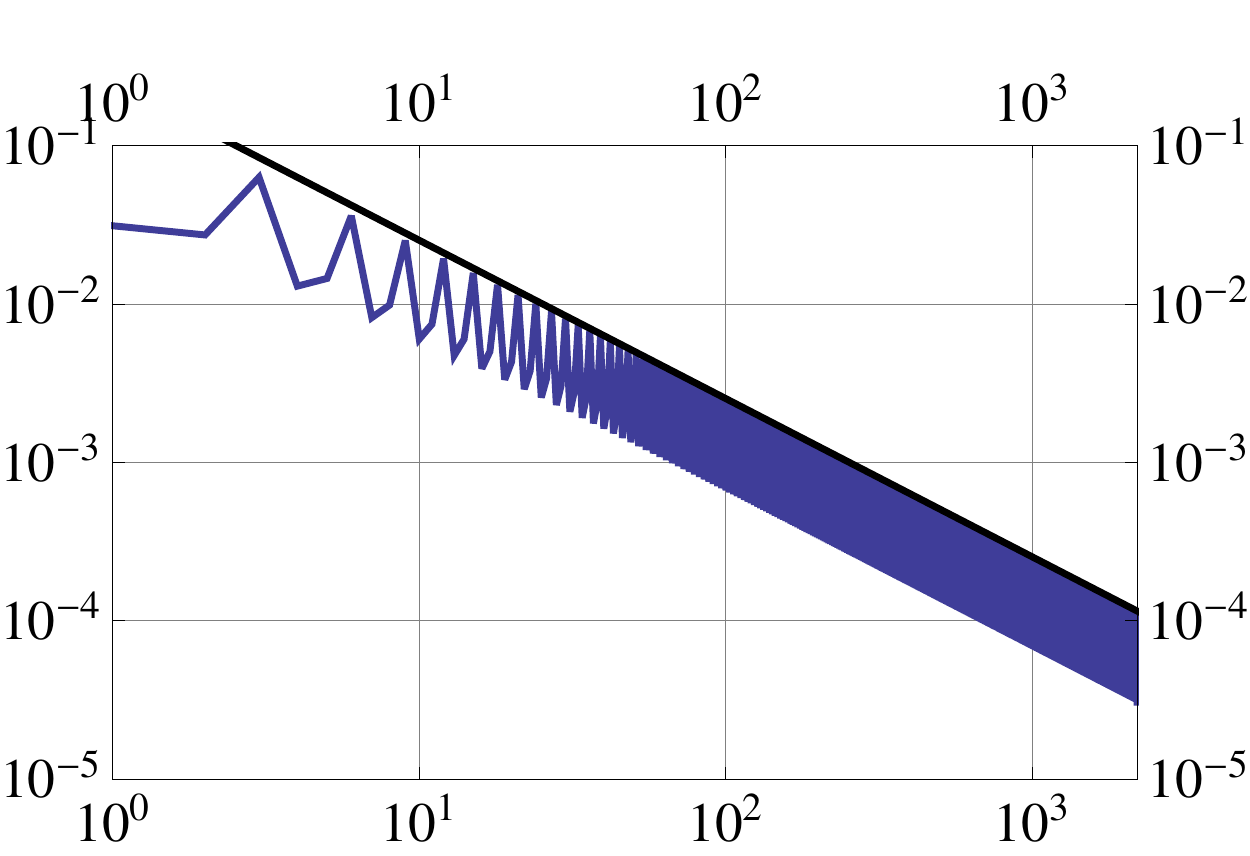}}
  \caption{$\varepsilon^{\prime }(x)$: Error at $x=a$, $\alpha (a)=1$ , $C(a)=0.25293$.
Absolute value vs polynomial order, $p=1,\ldots,2200$.}
  \label{fig:fig2}
\end{figure}
Figure~\ref{fig:fig2} shows the error for $x=a.$ We see $\alpha(1) =1$ and $C(a)=0.25293$.
The estimate (\ref{eq:4.9}) leads to $C(a)=0.31830$. Once more 
the rates agree but the observed constant is lower.

\begin{figure}[ht]
	\centering
	\subfloat[Error at $x=1/10$, $\alpha (1/10)=1$, $C(1/10)=0.84622$.]
	{\label{fig:fig3a}\includegraphics[width=0.45\textwidth]{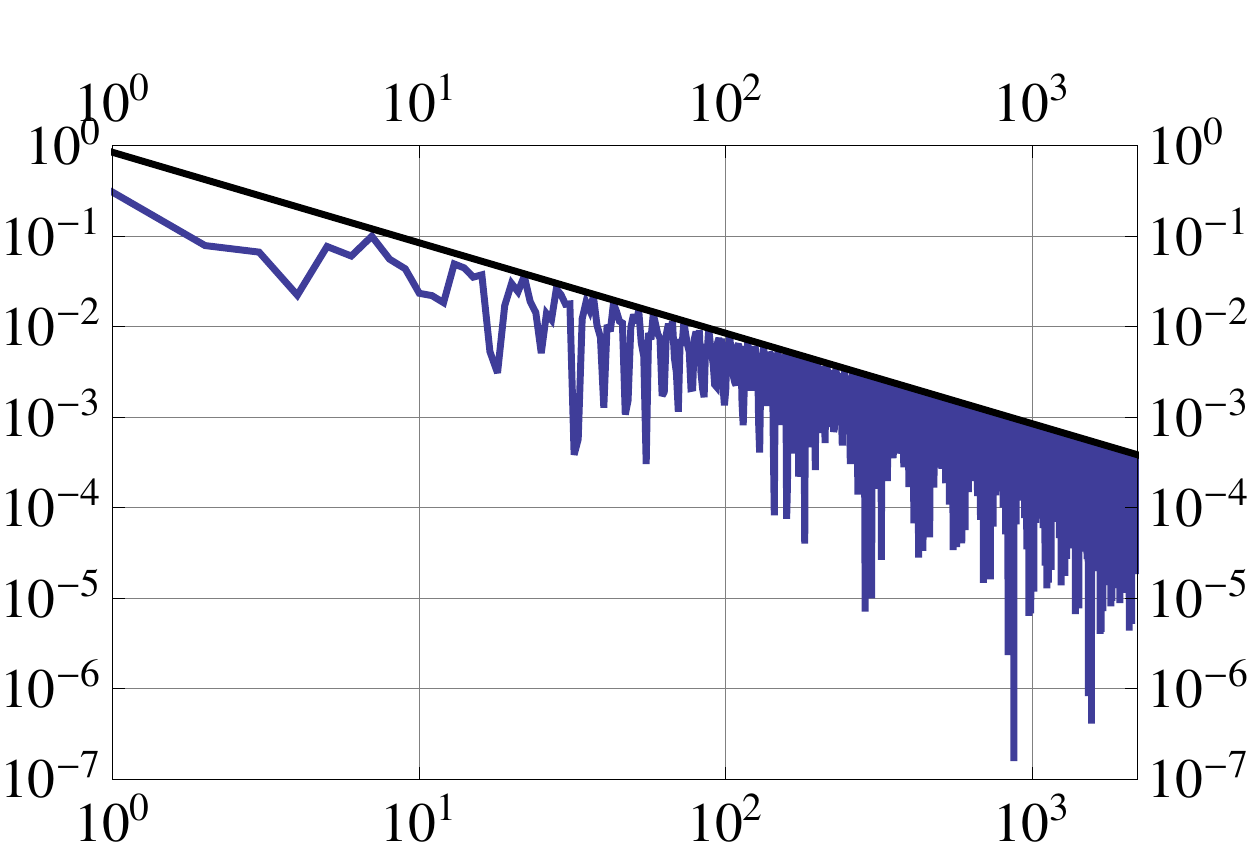}}\quad
	\subfloat[Error at $x_{2}=-1+10^{-2}$, $\alpha (x_{2})=1$, $C(x_{2})=0.889506$.]
	{\label{fig:fig3b}\includegraphics[width=0.45\textwidth]{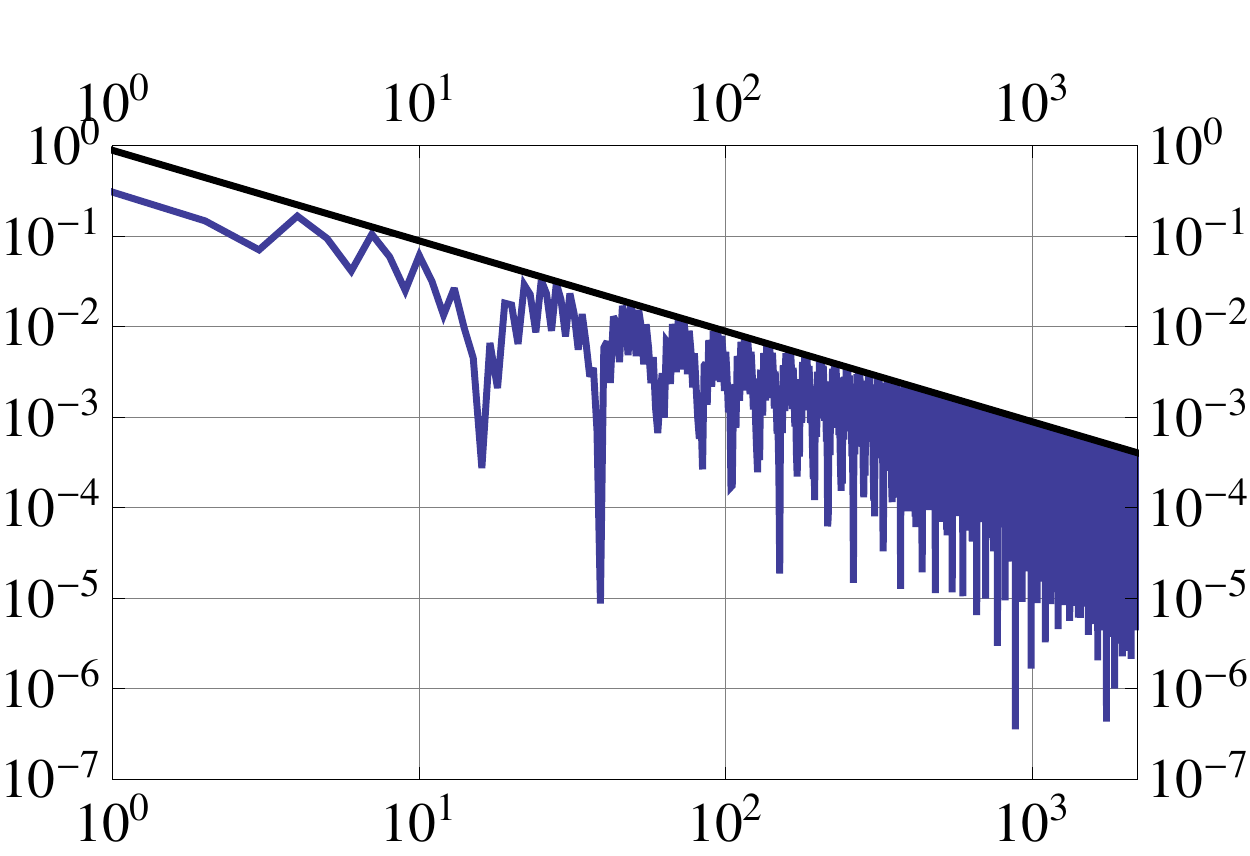}}\\
	\subfloat[Error at $x_{4}=-1+10^{-4}$, $\alpha (x_{4})=1$, $C(x_{4})=2.733292$.]
	{\label{fig:fig3c}\includegraphics[width=0.45\textwidth]{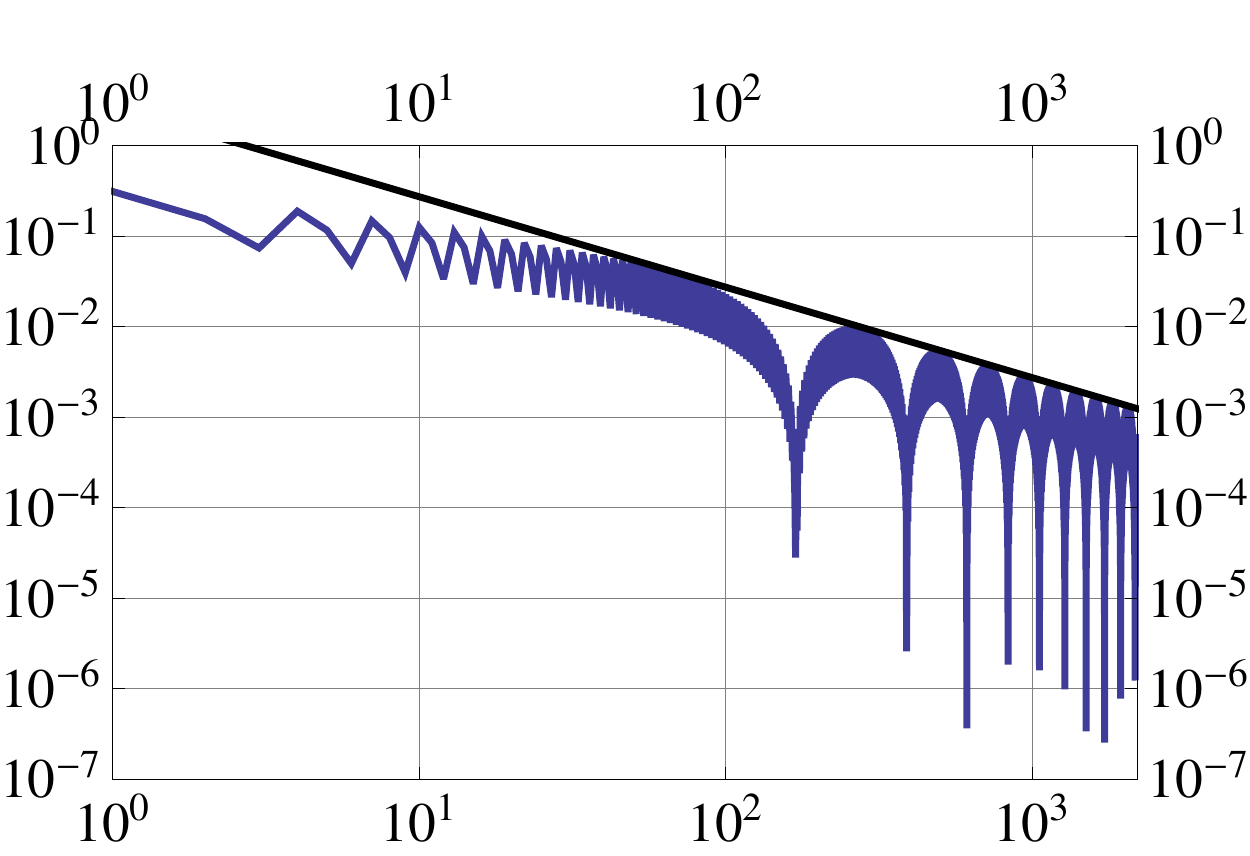}}\quad
	\subfloat[Error at $x_{6}=-1+10^{-6}$, $\alpha (x_{6})=1$, $C(x_{6})=7.765878$.]
	{\label{fig:fig3d}\includegraphics[width=0.45\textwidth]{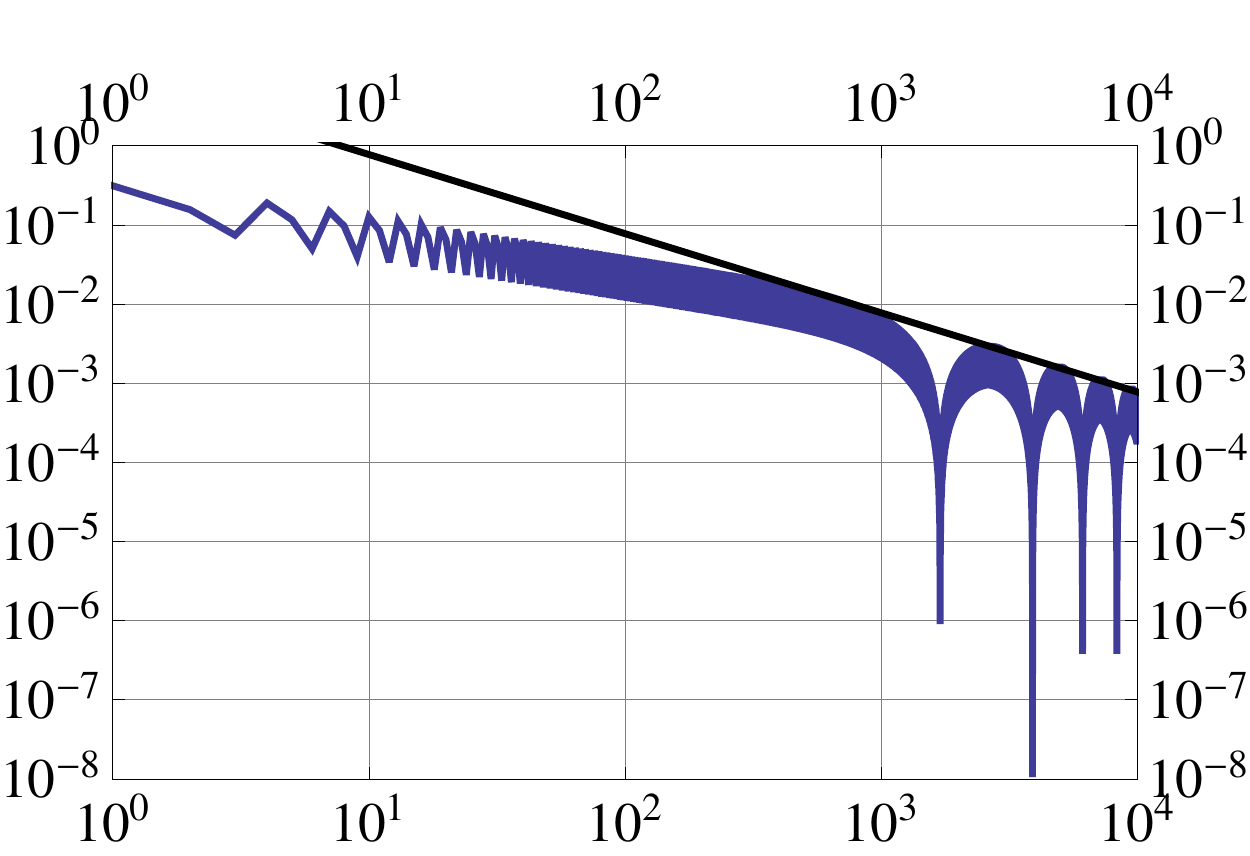}}
  \caption{$\varepsilon^{\prime }(x)$: Absolute value vs polynomial order, 
$p=1,\ldots,2200$, except for (d), $p=1,\ldots,10000$.}
  \label{fig:fig3}
\end{figure}
Figure~\ref{fig:fig3a} shows the error for $x=0.1$ with convergence rate  $\alpha (0.1)=1$
and $C(0.1)=0.84622$. This rate of convergence follows from the Theorem~\ref{thm:1}. 
By the Theorem~\ref{thm:1} we get $C=32.793$, i.e., largely overestimated error estimate
although with correct rate. 
Of course we have to have in mind that the
estimate in Theorem~\ref{thm:1} covers much larger class of functions and hence the
constants have to be larger.
Moreover, the overall convergence pattern has a very different character than previously
and the lower bound seems not have the same rate as the upper bound. 
We are not able to make any
hypothesis about the lower bound and there are no theoretical results available 
that are addressing the lower bound. 
\begin{figure}[ht]
	\centering
	\subfloat[Pattern at $x=-1$]
	{\label{fig:figB3a}\includegraphics[width=0.45\textwidth]{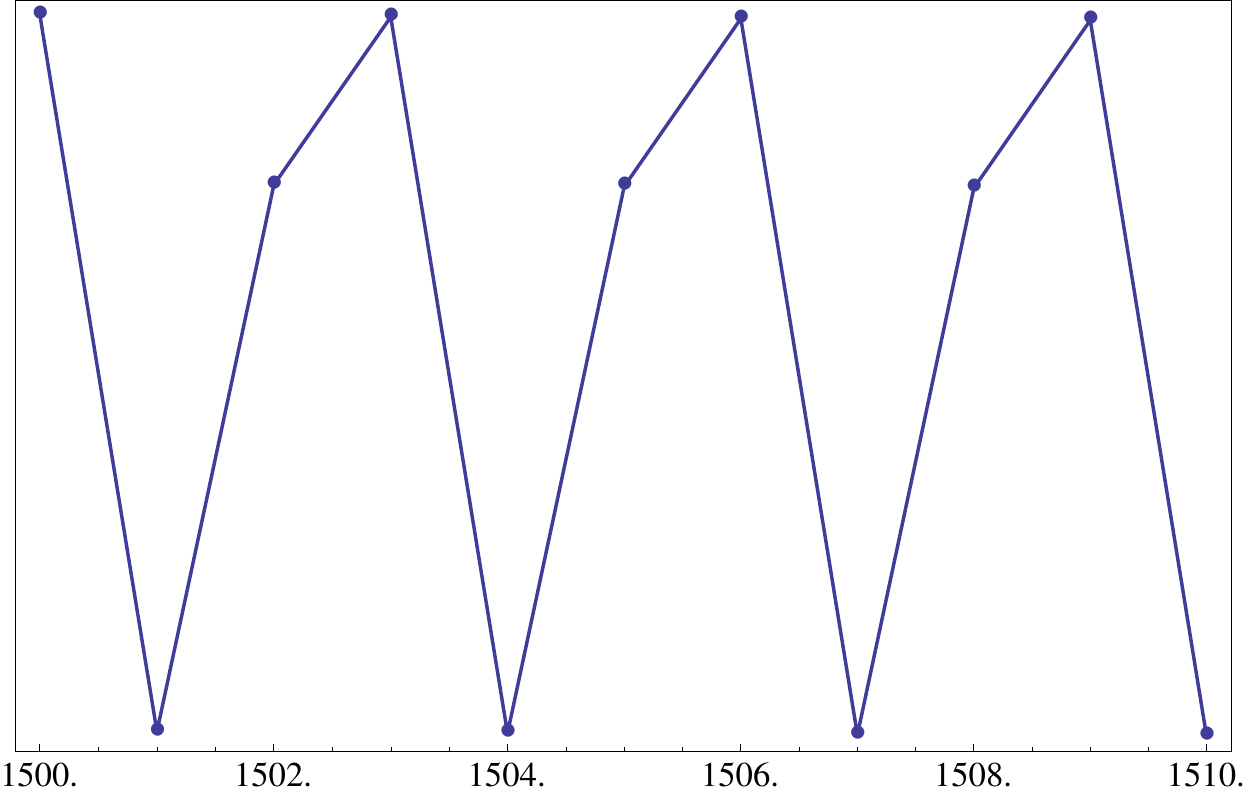}}\quad
	\subfloat[Pattern at $x=1/10$]
	{\label{fig:figB3b}\includegraphics[width=0.45\textwidth]{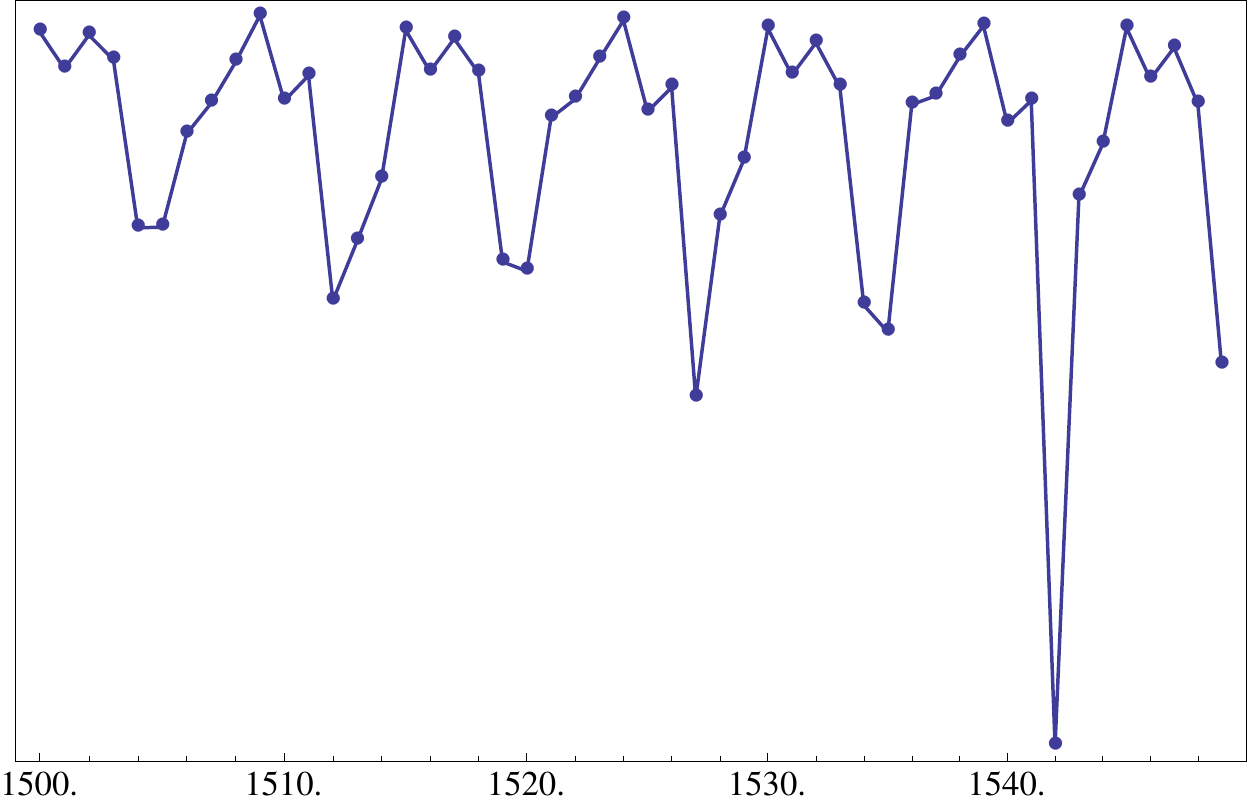}}
  \caption{Patterns: Absolute value vs polynomial order.}
  \label{fig:figB3}
\end{figure}

We see different rate of convergence for $x=0.1$ and $x=\pm 1.$
Hence let us
address the error  close to the boundary. Figures~\ref{fig:fig3b}--\ref{fig:fig3d} shows the error
for $x_{m}=-1+10^{-m},m=2,4,6$. We see that the pre-asymptotic range
increases with $m$ . We have rate $\alpha =1$ for all $x_{m}$ and $C$
increasing with $m$ $C(x_{2}) = 0.889506, C(x_{4}) =2.733292, C(x_{6}) =7.765878$. 
Figure~\ref{fig:figB3} shows very different pattern of the error $\mid
\varepsilon _{p}^{\prime }(x)\mid $ for $x=1$ and $x=0.1$ 
for $1000\leq p\leq 1500.$

\begin{figure}[ht]
	\centering
        \subfloat[Growth of the constant $C(-1+\xi)~\sim~\hat{D}\,\xi^{-1/4}$.]
        {\label{fig:figB4a}\includegraphics[width=0.45\textwidth]{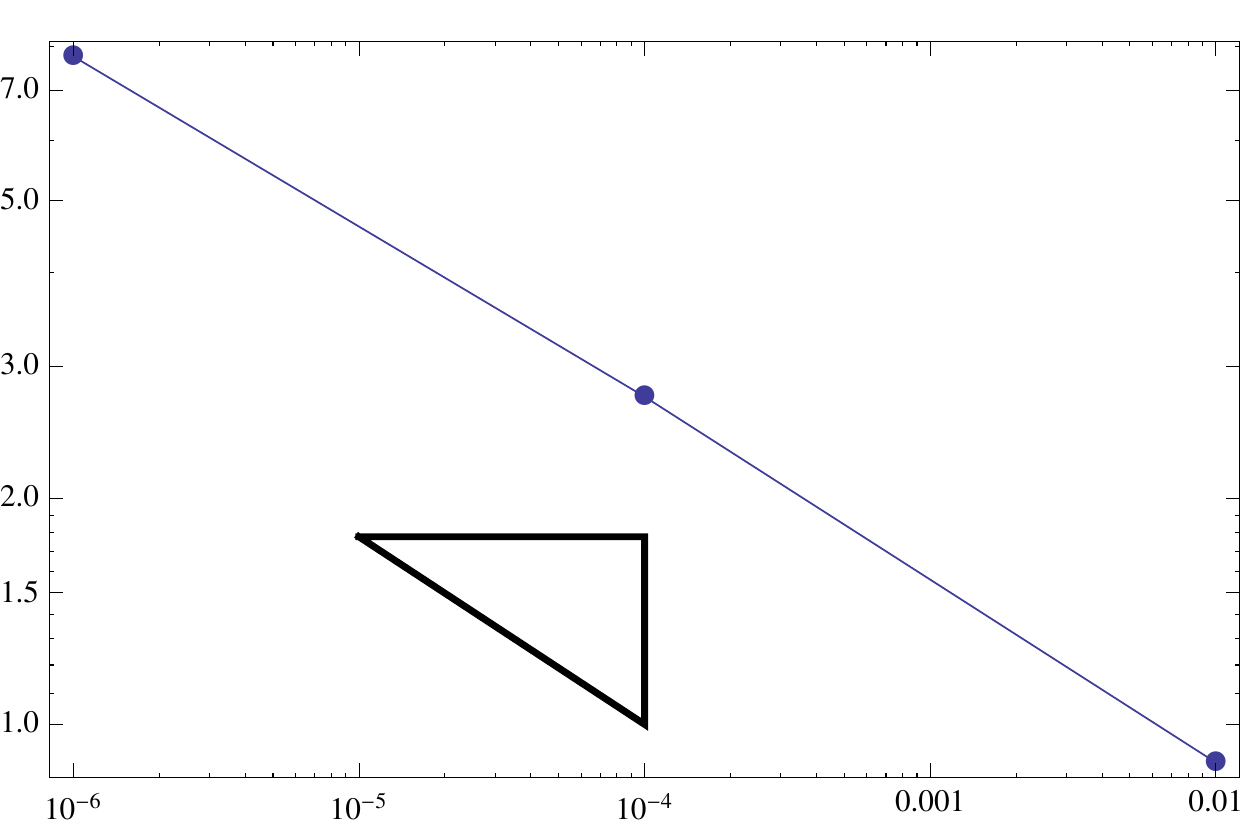}}\quad
        \subfloat[Growth of the constant $C(a+\xi)~\sim~\hat{D}\,\xi^{-1}$.]
	{\label{fig:figB4b}\includegraphics[width=0.45\textwidth]{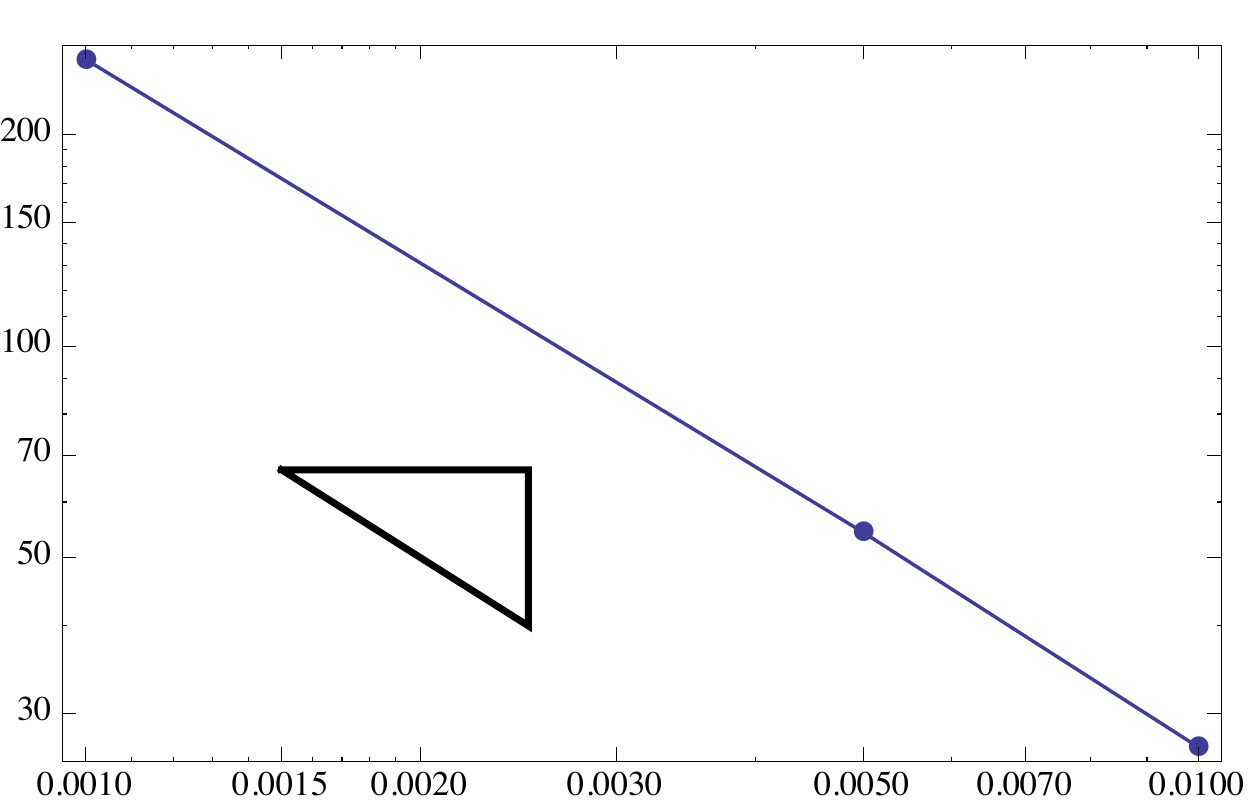}}
  \caption{$\varepsilon^{\prime }(x)$: Value of the coefficient vs distance to the point of interest.}
  \label{fig:figB4}
\end{figure}

The growth of $C$ is caused by the different rates $%
\alpha (-1)<\alpha (x_{m})$ for all $m$. In the Theorem~\ref{thm:1} we have seen the term 
$(1-x^{2})^{-3/2}$ which indicates the growth of the rate.  In the
Figure~\ref{fig:figB4a} we show in the log log scale the growth of the constant  $C(x)$. 
We see that 
$C(-1+\xi)\sim D(-1+\xi)^{\beta }$ with $\beta =-1/4$ and small $\xi > 0$.
From the theory we have $\beta = -3/2$.

\begin{figure}[ht]
	\centering
	\subfloat[Detail of $u^{\prime }(x)$ at $p=1000$.]
	{\includegraphics[width=0.45\textwidth]{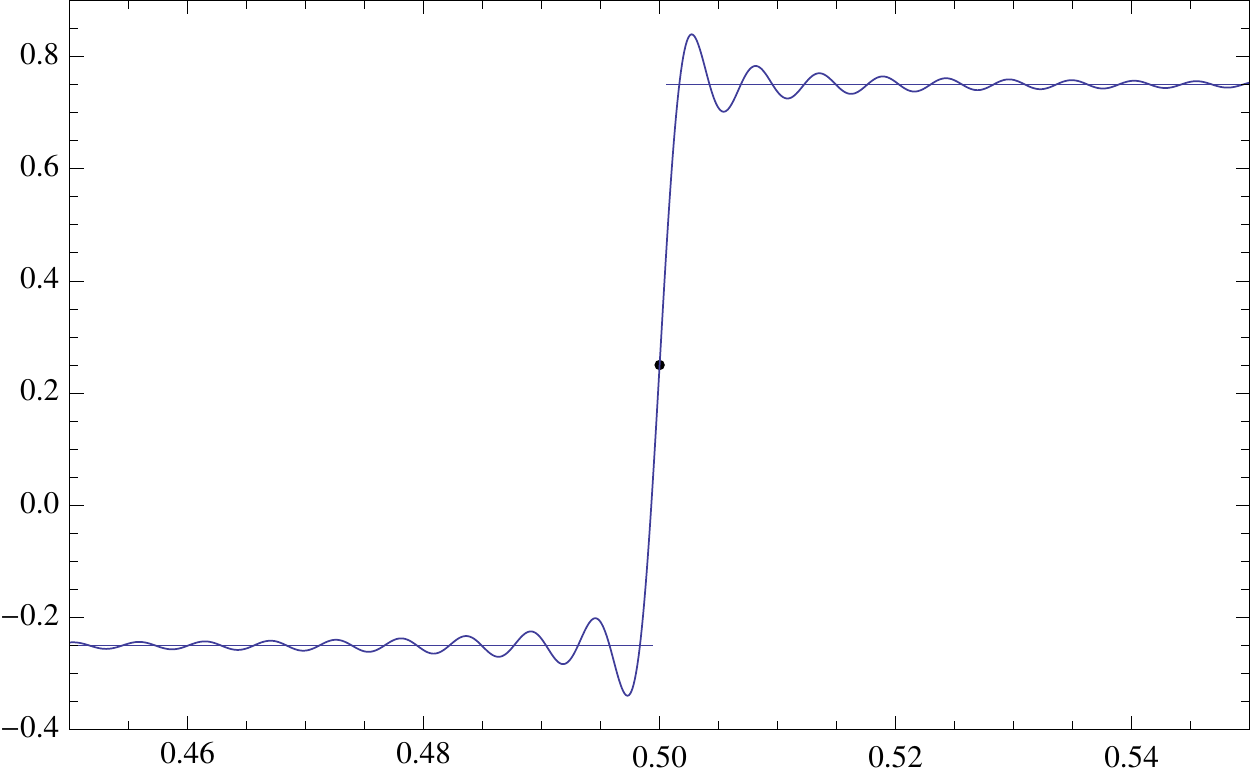}}\quad
	\subfloat[Detail of $u^{\prime }(x)$ at $p=2000$.]
	{\includegraphics[width=0.45\textwidth]{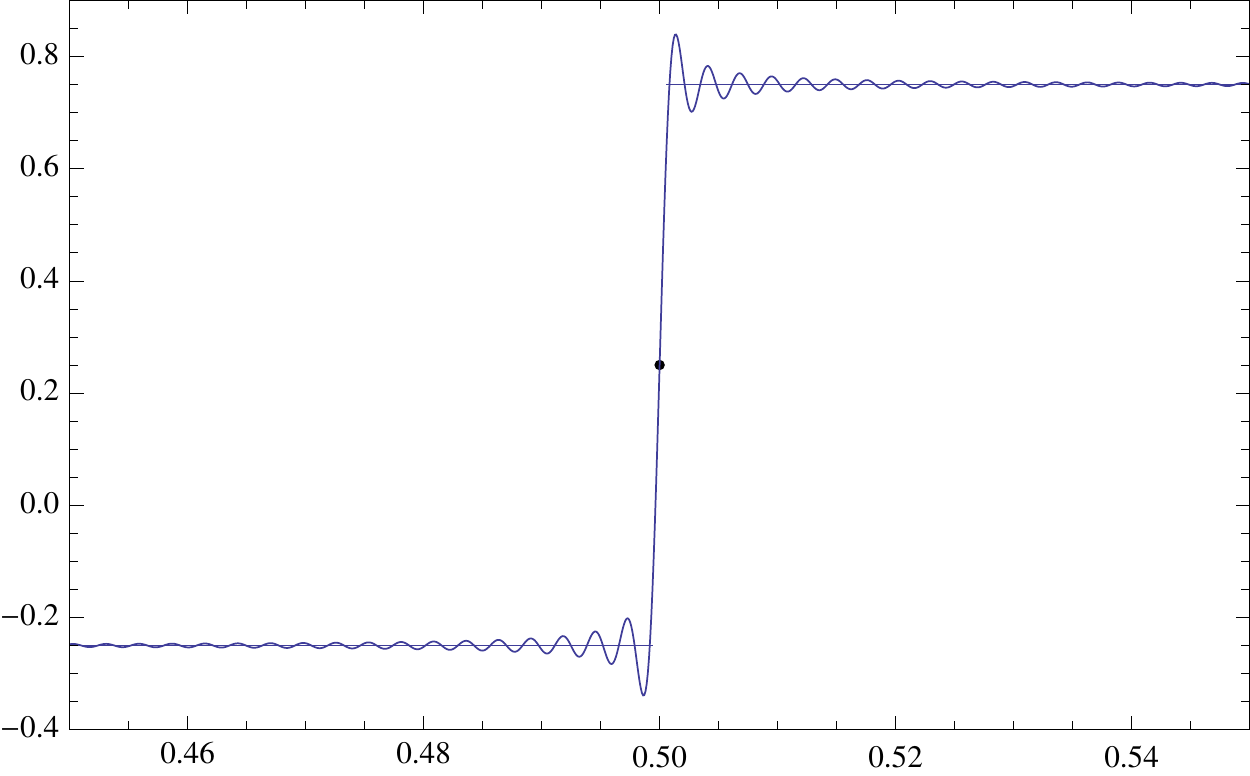}}
  \caption{Gibbs phenomenon}
  \label{fig:fig4}
\end{figure}
In the point $x=a$ the function $u'$ is discontinuous. Nevertheless
convergence to $\frac{1}{2}(u'(a+0)+u'(a-0))$ and elsewhere in the
neighborhood of $x=a$ has the same rate namely $\alpha =1$. In the Figure~\ref{fig:fig4} 
we see a typical error overshoot which is independent of $p$.
Denoting by $y$ the position of the maximal error we have 
$\mid y-a\mid =D\,p_{{}}^{-1}$ with $D=2.7777$.
This is the well known Gibbs phenomenon.
We have $\varepsilon ^{\prime }(a+\xi )\leq D\,\xi ^{-1}p^{-1}$
with $D$ independent of $p$ and small $\xi$.
\begin{figure}[ht]
	\centering
	{\includegraphics[width=0.45\textwidth]{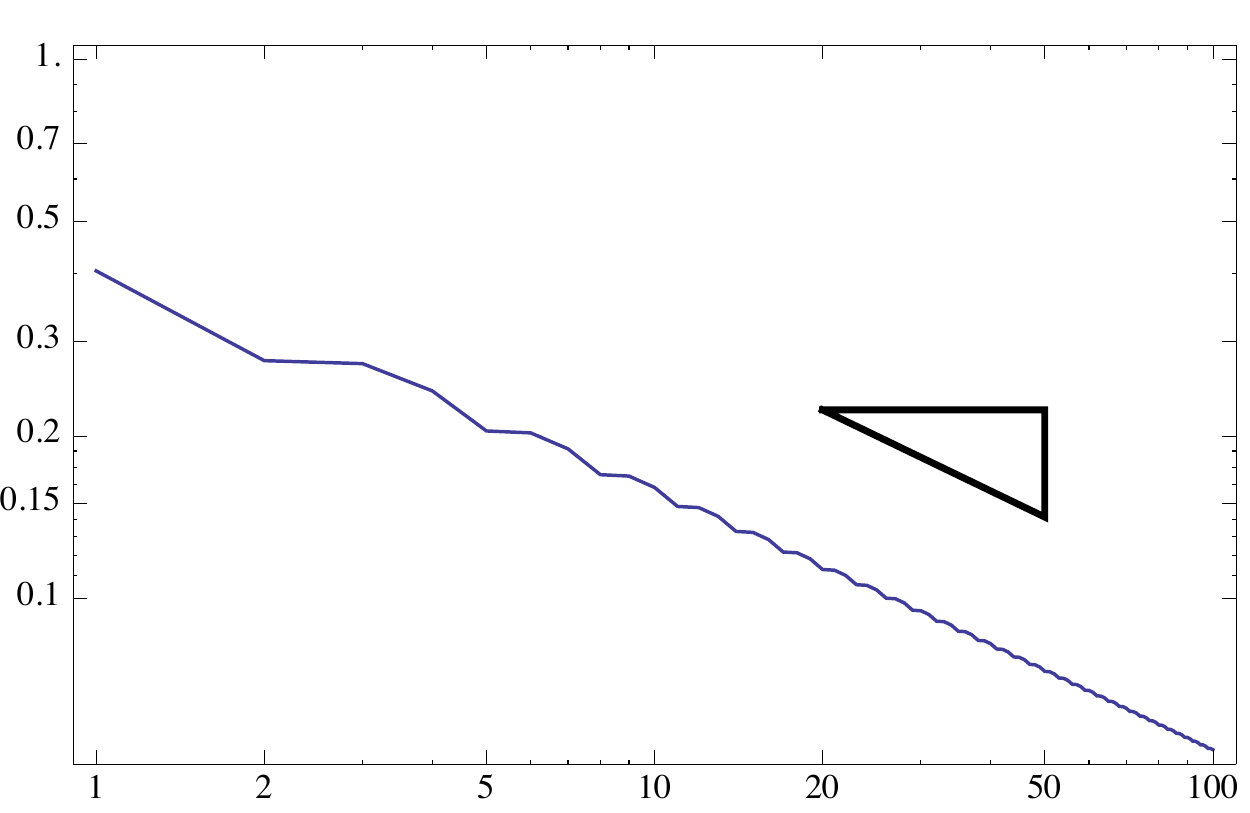}}
  \caption{Convergence in the norm: $\parallel \varepsilon _{p}\parallel _{\mathcal{E}}\leq C\,p^{-1/2}$. Value vs polynomial order, $p=1,\ldots,100$.}
  \label{fig:fig5}
\end{figure}
In the Figure~\ref{fig:fig5} we show the convergence of $\parallel \varepsilon _{p}^{\prime
}\parallel _{L^{2}}=\parallel \varepsilon _{p}\parallel _{\mathcal{E}}$.
As expected, since $u^{\prime }\in
B_{2,\infty }^{1/2}$ it follows from the theory that 
$\parallel \varepsilon _{p}\parallel _{\mathcal{E}}\leq
Cp^{-1/2}$ because of the regularity of $u$. This
coincides very well with the numerical computations.

\begin{remark}
Approximation of the function $u^{\prime }(x)$ was
  first addressed already more than 100 years ago \cite{c14}. 
  It has been used as an example for the Legendre
expansion in the most simple setting in various books. Let us mention
for example \cite[p. 163]{c9} and \cite[p. 58]{c9a}.
\end{remark}
Let us summarize our results
\begin{itemize}
  \item 
    [1.]{The classical error estimate in the energy norm is in a very good
agreement with the numerical results }

  \item 
    [2.]{The rate of convergence for }$x\in I,x\neq -1,1,a$, {does
follow from the general theory based on only the total variation of the
function but the constant is very inaccurate.} The classical Gibbs phenomenon 
is clearly visible in the Figure~\ref{fig:fig4}.

We see that the convergence of the Legendre polynomials in a point 
$x$ is not governed by the smoothness of the function in
that point and its neighborhood. In our case the function was constant in
the neighborhood. There is a strong pollution effect. 
We note that this pollution could be removed by a postprocessing \cite{c13}.
It is
characteristic for the $p$-version that the pollution in the
boundary points is larger than in its neighbors which leads to the boundary
layer in the convergence. This is a significant difference in comparison to the $h$-version. 
For the analysis of the pollution in the $h$-version we refer to \cite[Sect 9]{c7}.

  \item 
    [3.]We have seen that in the neighborhood of the points $x=\pm 1$ 
and $x=a$ the error $\varepsilon_p ^{\prime }(x)$
behaves differently. This behavior can be described by using weighted space
$L_{w}^{\infty }$ with the norm $\parallel \varepsilon _{p}^{\prime
}\parallel _{L_{w}^{\infty }}=\max_{x\in (-1,1)}\mid \varepsilon_p ^{\prime
}(x)w(x)\mid $ with $w(x)=\mid 1-x\mid ^{\alpha }\mid 1+x\mid ^{\beta }\mid
x-a\mid ^{\gamma }$. Particularly in our case we have $\parallel
\varepsilon _{p}^{\prime }\parallel _{L_{w}^{\infty }}\leq Cp^{-1}$ 
with $\alpha =\beta =1/2$ and $\gamma =1$. Nevertheless
this characterization gives no information about the behavior in the
singular points. 
\end{itemize}

\section{Legendre expansion of the solution $u$ given in (\ref{eq:2.2})}
\label{sec:legendreu}

Let $u(x)=\sum_{k=0}^{\infty }c_{k}P_{k}(x)$ 
be the Legendre expansion of the solution $u(x)$ of (2.1). 
Then $\sum_{k=0}^{p}c_{k}P_{k}$ $(x)$ is not the $p$-version
approximate solution of $u(x)$ because the constraint 
$u(\pm 1)=0$ would be not satisfied. To prevent any
misunderstanding we will write $w(x)$ instead of $u(x)$ and 
$w_{p}(x)=\sum_{k=1}^{p}c_{k}$ $P_{k}(x)$ with $\eta
_{p}(x)=w(x)-w_{p}(x)$ the error. Obviously 
\[
  u(x)=w(x)=\mid x-a\mid +\text{ linear function}
\]
and hence $\eta _{p}(x)$ is the error
of the Legendre expansion of the function $\mid x-a\mid $. The error
of the $p$-version will be analyzed in the next section.

It is easy to see that the
coefficients of the expansion of $u(x)$ are with $a_k$ given in (\ref{eq:4.3})
\begin{equation}
c_{0}=-a_{1}/3, c_{k}=-a_{k+1}/(2k+3)+a_{k-1}/(2k-1), k=1,2,3,\ldots
\end{equation}

\begin{figure}[ht]
	\centering
	\subfloat[Error at $x=-1$, $\alpha (-1)=3/2$, $C(-1)=0.42625$.]
	{\includegraphics[width=0.45\textwidth]{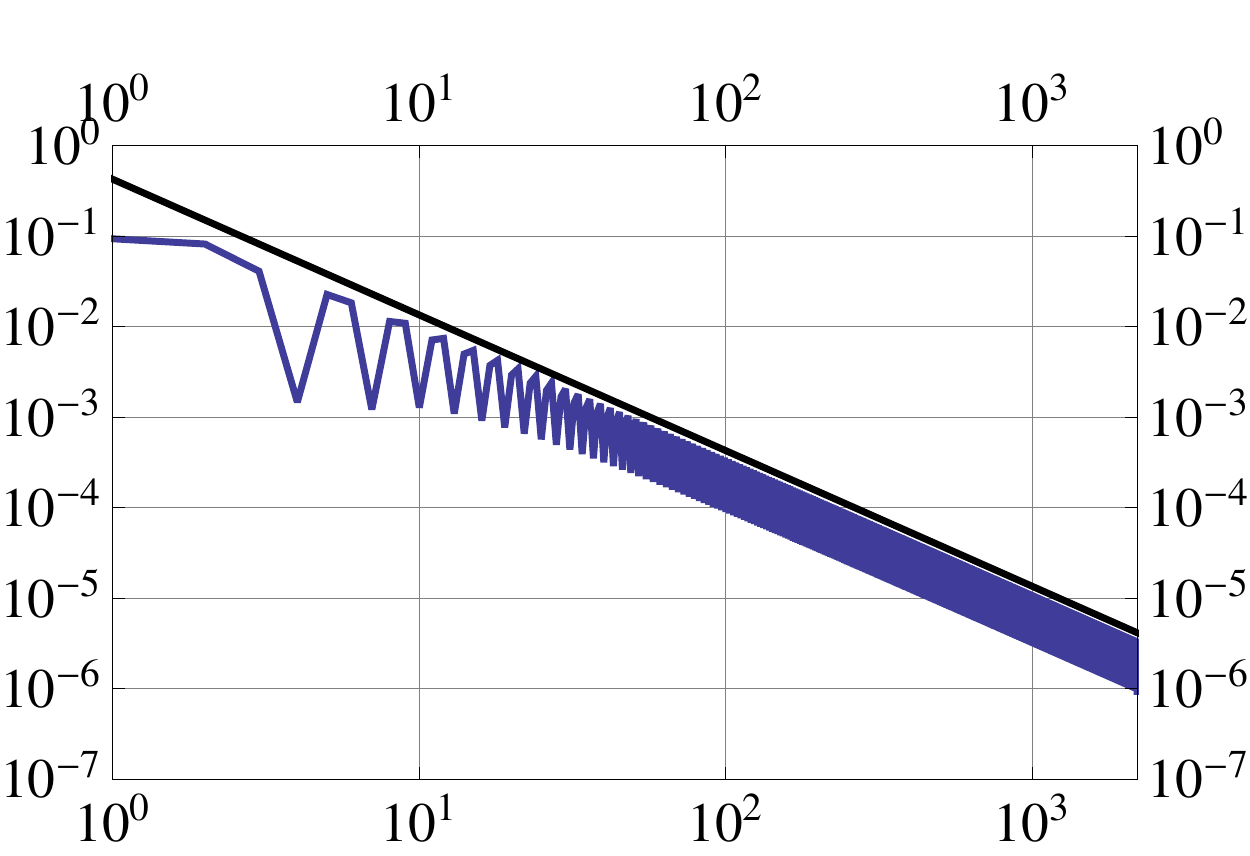}}\quad
	\subfloat[Error at $x=-99/100$, $\alpha (-99/100)=2$, $C(-99/100)=0.76483$.]
	{\includegraphics[width=0.45\textwidth]{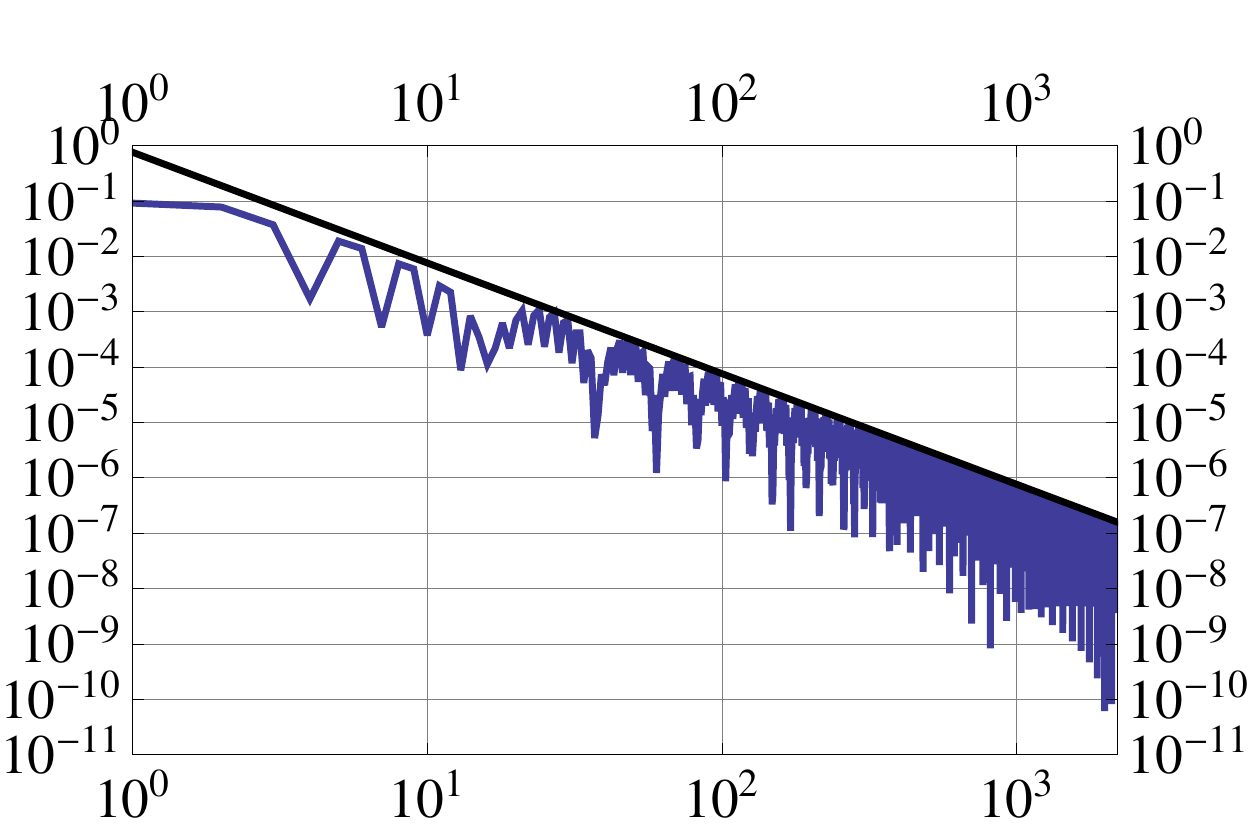}}\\
	\subfloat[Error at $x=1/10$, $\alpha (1/10)=2$, $C(1/10)=0.73185$.]
	{\includegraphics[width=0.45\textwidth]{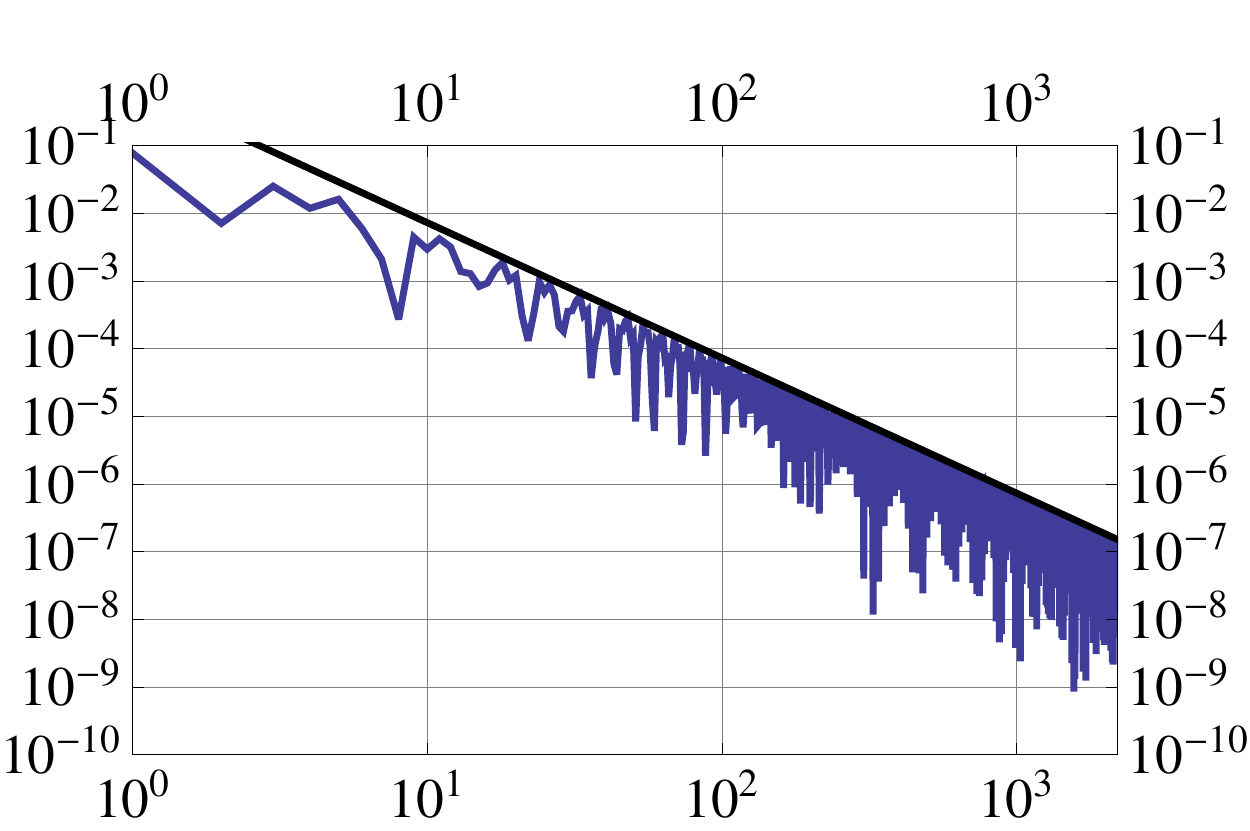}}
  \caption{$\eta_{p}(x)$: Absolute value vs polynomial order, $p=1,\ldots,2200$.}
  \label{fig:fig6}
\end{figure}

\begin{figure}[ht]
	\centering
	\subfloat[Error at $x=a$, $\alpha (a)=1$, $C(a)=0.274738$. $p_{\max}=10000$.]
	{\label{fig:fig7a}
\includegraphics[width=0.45\textwidth]{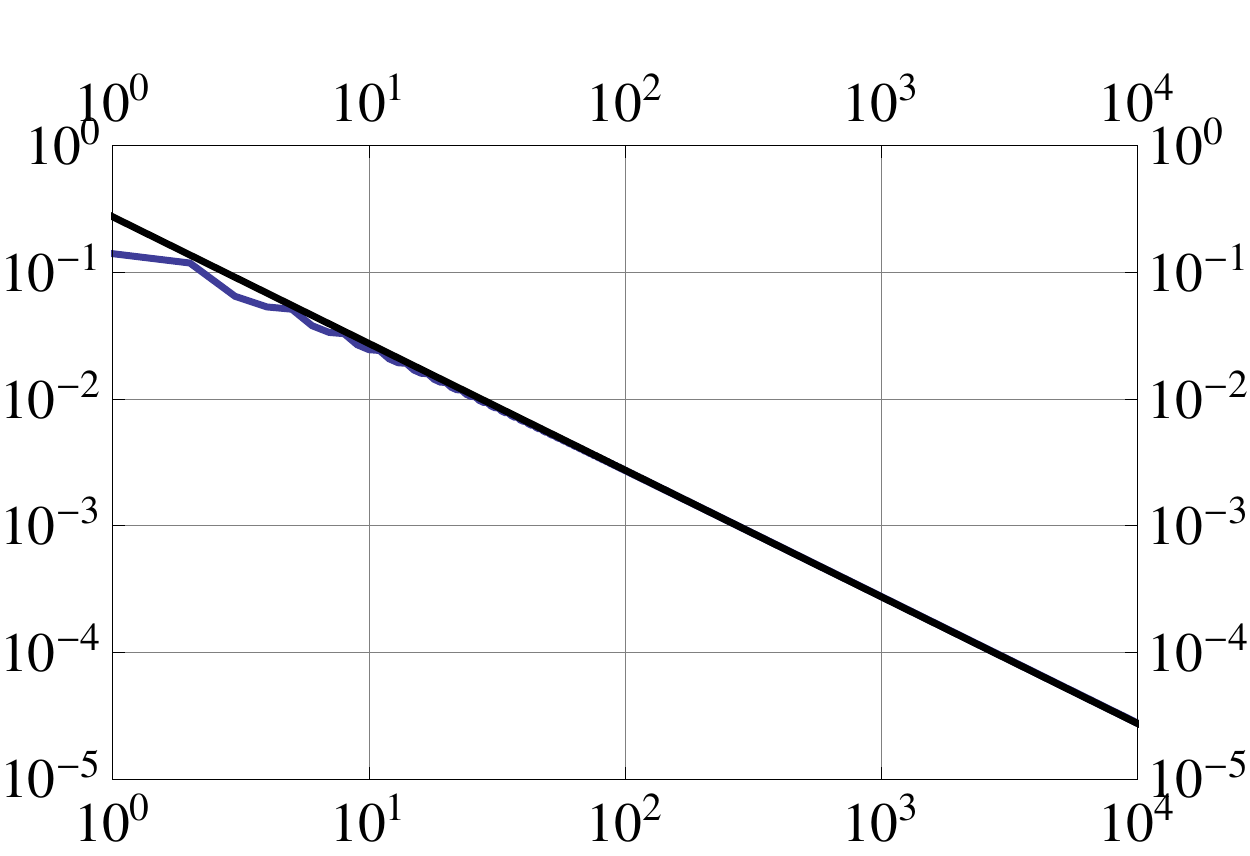}}\quad
	\subfloat[Error at $x=a+1/100$, $\alpha (a+1/100)=2$, $C(a+1/100)=24.5325$. $p_{\max}=2200$.]
	{\includegraphics[width=0.45\textwidth]{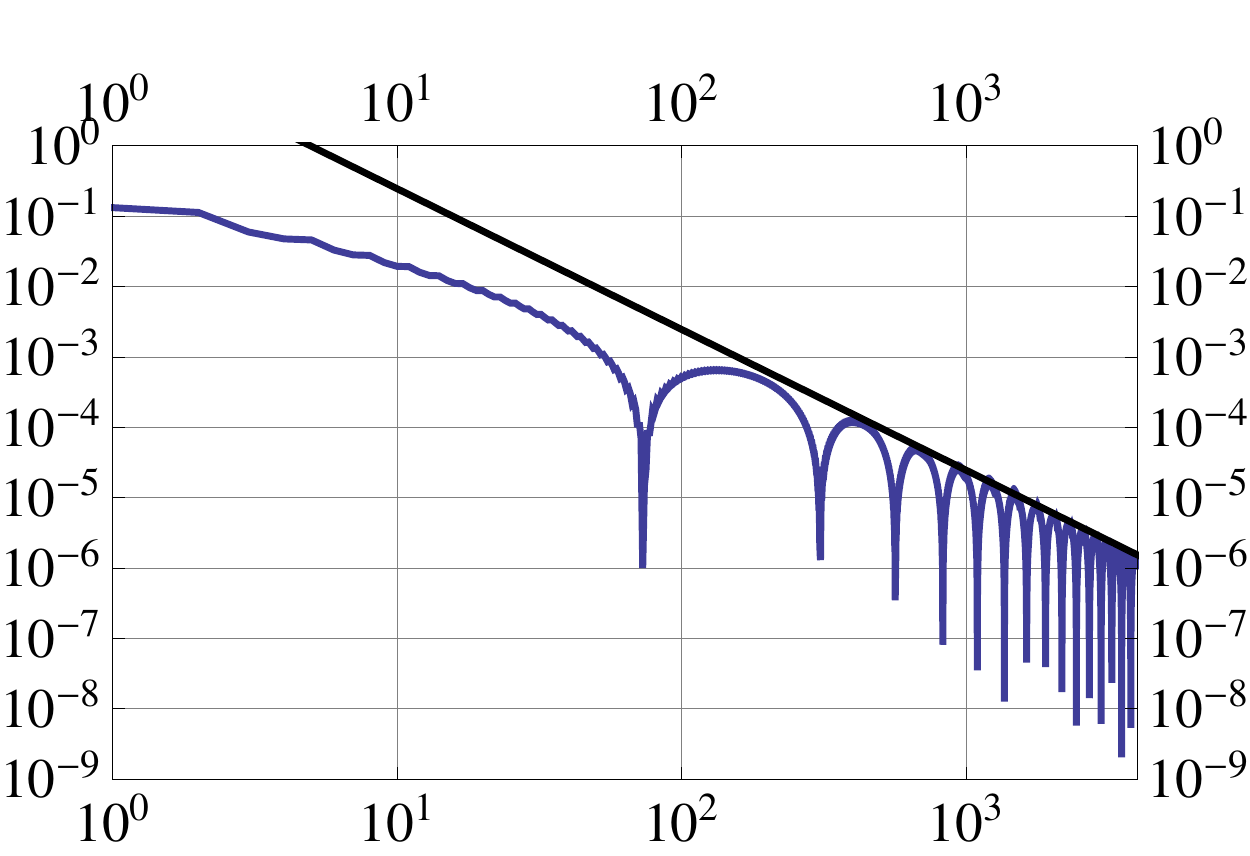}}
        \caption{$\eta_{p}(x)$: Absolute value vs polynomial order, $p=1,\ldots,p_{\max}$.}
  \label{fig:fig7}
\end{figure}
\begin{figure}[ht]
	\centering
        \subfloat[Growth of the constant $C(-1+\xi)~\sim~\hat{D}\,\xi^{-1/4}$.]
	{\includegraphics[width=0.45\textwidth]{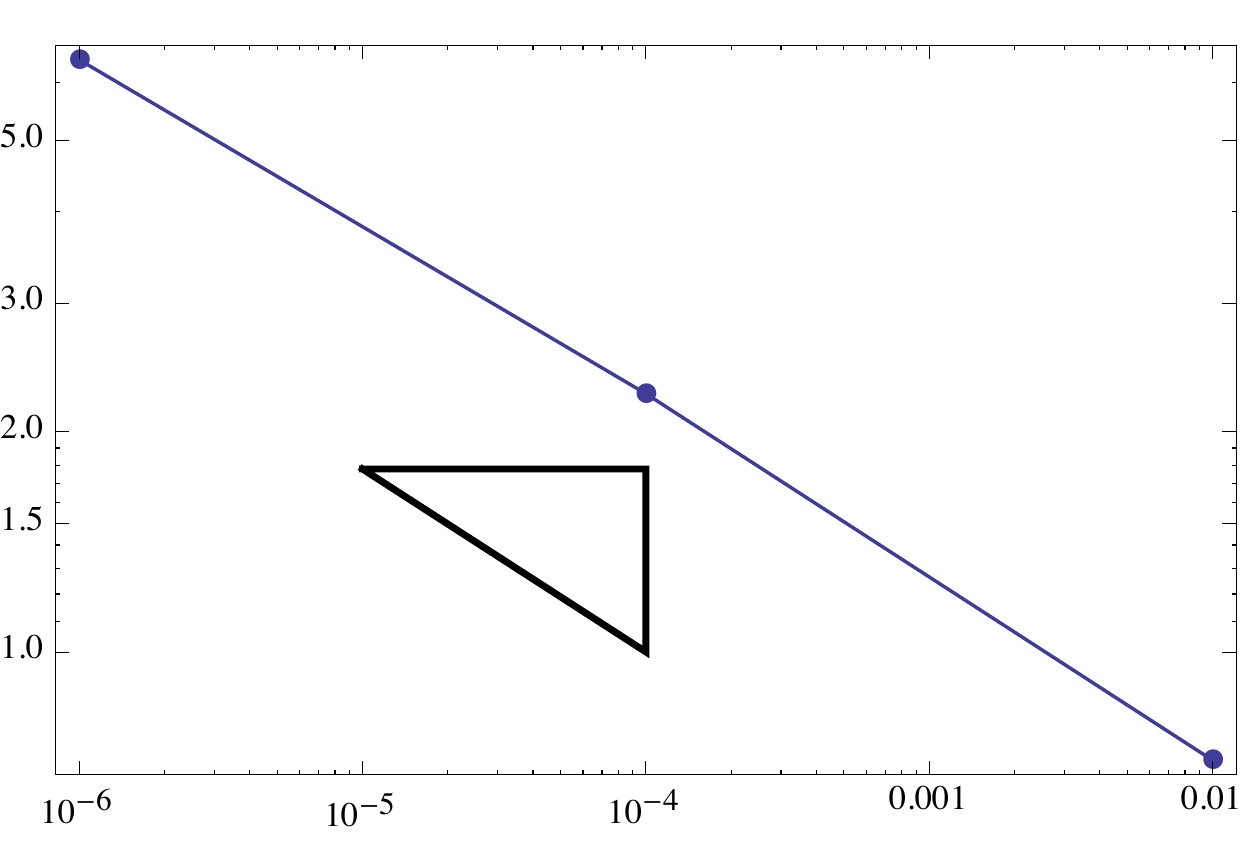}}\quad
         \subfloat[Growth of the constant $C(a+\xi)~\sim~\hat{D}\,\xi^{-1}$.]
	{\includegraphics[width=0.45\textwidth]{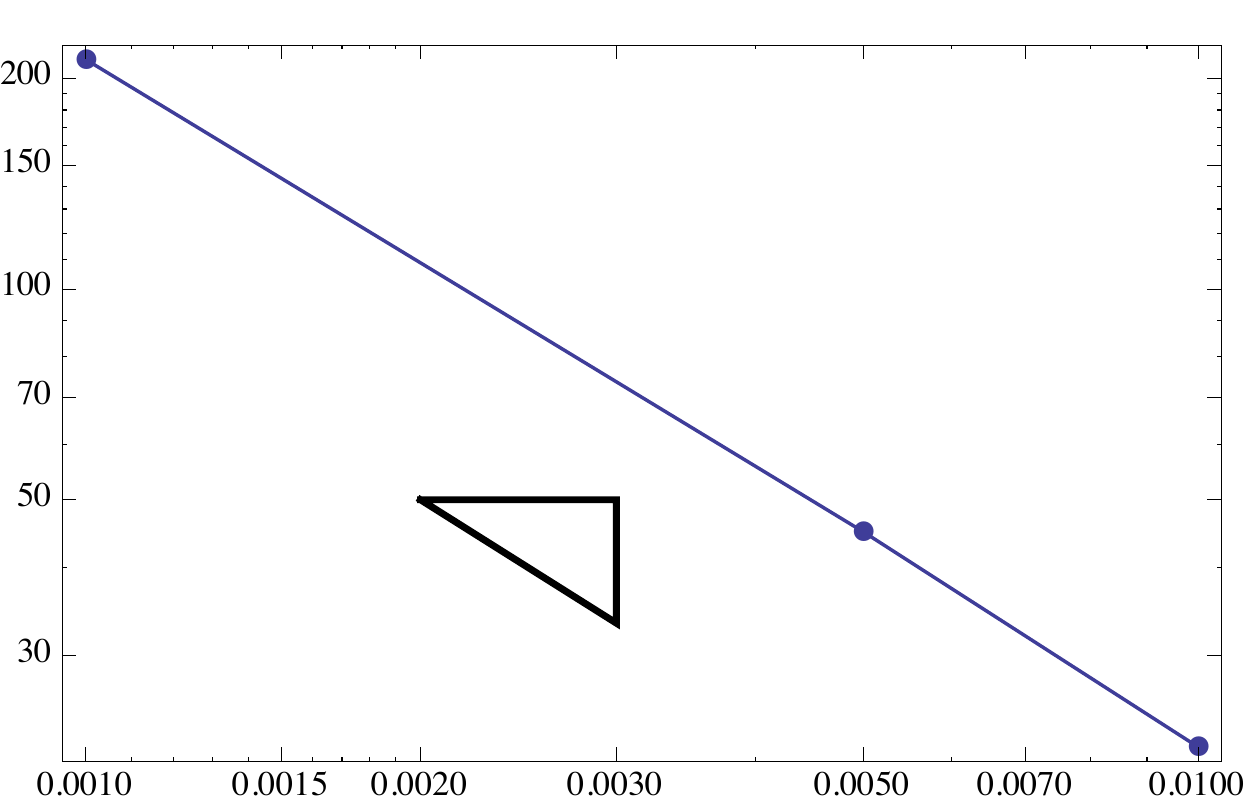}}
        \caption{$\eta_{p}(x)$: Value of the coefficient vs distance to the point of interest.}
  \label{fig:figB7}
\end{figure}
The function $u$ is smoother than $u'$ of the previous section. 
Numerically we see some analogous behavior
with $\alpha (x)=2$ for $x\neq a,\pm 1,$ $\alpha (\pm 1)=\frac{3}{2}$ 
and $\alpha (a)=1$. Once more we see the different rates of convergence for 
$x\neq a,\pm 1$. At $x=\pm 1$ the rate is by $\frac{1}{2}$ smaller then in
its neighboring points and we have similar increase of $C$ (Figure~\ref{fig:figB7}). 
In fact, the observed growth rate is exactly the same, $\beta = -1/4$, yet as expected,
the value of the constant is smaller here.
Figure~\ref{fig:fig6} shows
the $\mid \eta _{p}(x)\mid $ for $x=-1,0.99,0.1$. We have no Gibbs
phenomenon in its classical form but we have different rates of convergence
in $x=a$ and its neighboring points. The difference of the rates is exactly 1. This
difference in the rates is stronger here than at the boundary points. Figure~\ref{fig:fig7}
shows the error for $x=a$ and $x=a+0.01$ and the Gibbs phenomenon is
different. We have $\eta (a+\xi )\sim \xi^{-1} p^{-2}$.

Application of the Theorem~\ref{thm:1} leads to $\alpha (x)=1$ for $x\neq \pm 1,a$. 
For $x=a$ we get the error estimate $\mid \eta _{p}(a)\mid \leq C\,p^{-1}\lg p$,
the right rate up to the log term. Notice, that in Figure~\ref{fig:fig7a} the series
has been evaluated upto $p=10000$ and there is no evidence of the log term affecting
the convergence.

Typical theorem related to $\eta_{p}$ is
\begin{theorem}[\cite{c9}]
Let $f(x)$ satisfy on [-1,1] the Lipschitz \ condition with $\gamma >1/2$.
Then we have
\begin{equation}
\mid f(x)-\sum_{k=0}^{p}c_{k}P_{k}(x)\mid \leq
	c(x)/p^{\gamma -1/2}
\end{equation}
\end{theorem}
In our case we have $\gamma =1$ and the Theorem~\ref{thm:1} predicts the rate $\alpha
=1/2$ while we have seen $\alpha =1$.

\begin{figure}[ht]
	\centering
	{\includegraphics[width=0.45\textwidth]{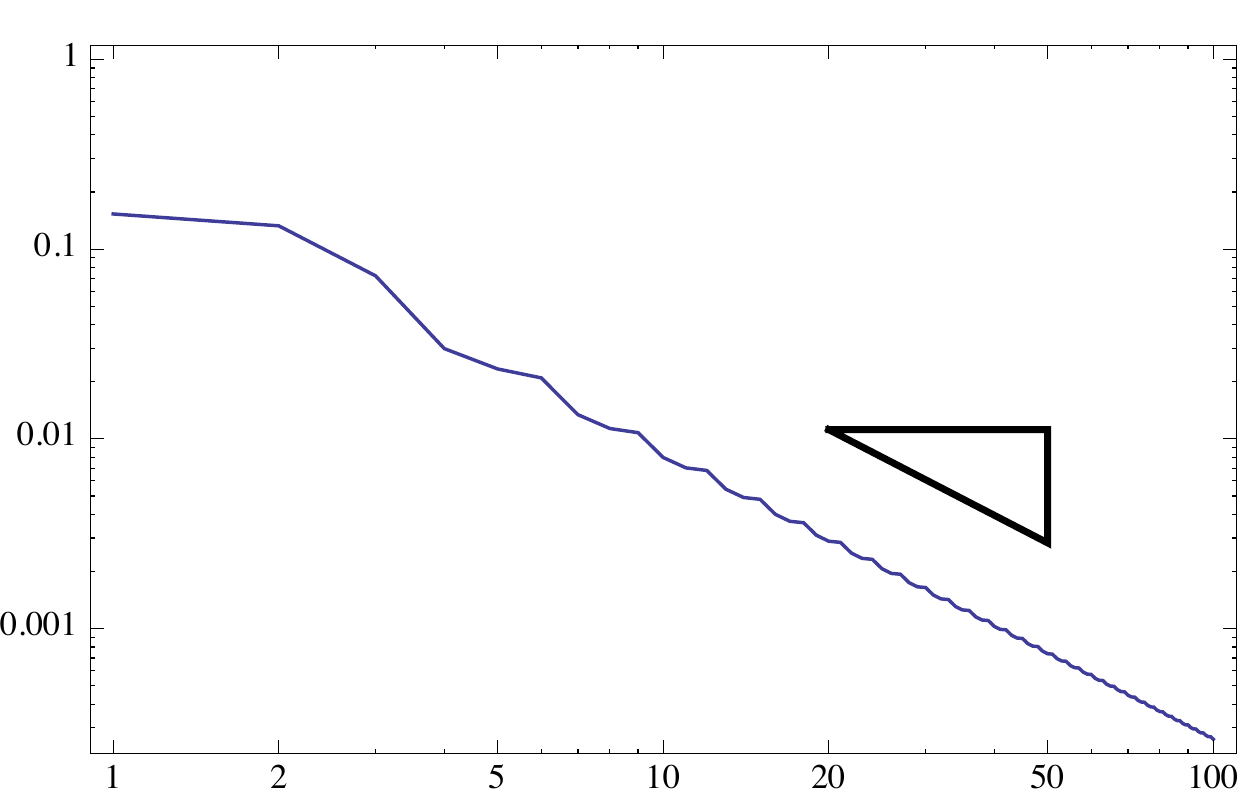}}
  \caption{Convergence in the norm: $\parallel \eta _{p}\parallel
_{L^{2}}\leq  Cp^{-3/2}$. Value vs polynomial order, $p=1,\ldots,100$.}
  \label{fig:fig8}
\end{figure}
In the Figure~\ref{fig:fig8} we show the convergence rate $\parallel \eta _{p}\parallel
_{L^{2}}\leq C\,p^{-3/2}$ as expected because of the regularity of the
function $u$.

Let us summarize the results.
\begin{itemize}
  \item[1.] The results are very similar as in the previous case -- only the
rate is increased by one except for $x=1$ because we addressed the
convergence to the average $(u'(a+0)+u'(a-0))/2$. This increase was caused by the increase 
of the smoothness of the expanded function.

\item[2.] We have seen that the theory gives more pessimistic results than
observed. The reason is that the theory deals with the functions having only
bounded variation. In addition the above theorems are
addressing the $L^{\infty }$ norm which does not distinguish
between the interior points and the points at the boundary.
\end{itemize}

\section{The error $\varepsilon_{p}$ of the $p$-version}
We addressed in the Section 5 the error of the
partial Legendre expansion $w_{p}(x)$. We underlined that $w_{p}(x)$
is not the $p$-version solution of the problem (2.1) because the
constraint of the boundary condition was not used.

The $p$-version solution $u_{p}$ is a modification of $w_{p}$. We get
\begin{equation}
\mid u_{p}(x)\mid =\sum_{k=0}^{p+1}b_{k}P_{k}(x)
\end{equation}
where
\begin{equation}
\begin{array}{ll}
b_{k}=a_{k-1}/(2k-1), &\mbox{ (instead of $c_{k}=a_{k-1}/(2k-1)-a_{k+1}/(2k+3))$}, \\
b_{1}$ $=-a_{2}/5, &\mbox{ (instead of $(c_{1}=a_{1}/3-a_{3}/7)$)}, \\
b_{0}=-a_{1}/3, &\mbox{ (unchanged)}.
\end{array}
\end{equation}
Above in the parentheses we list the coefficients $c_{k}$ of the direct
expansion of $w$ addressed in the previous section.

\begin{figure}[ht]
	\centering
        \subfloat[Error at $x=-1+10^{-6}$, $\alpha(-1+10^{-6})$=2, $C(-1+10^{-6}$)=0.0013.]
	{\includegraphics[width=0.45\textwidth]{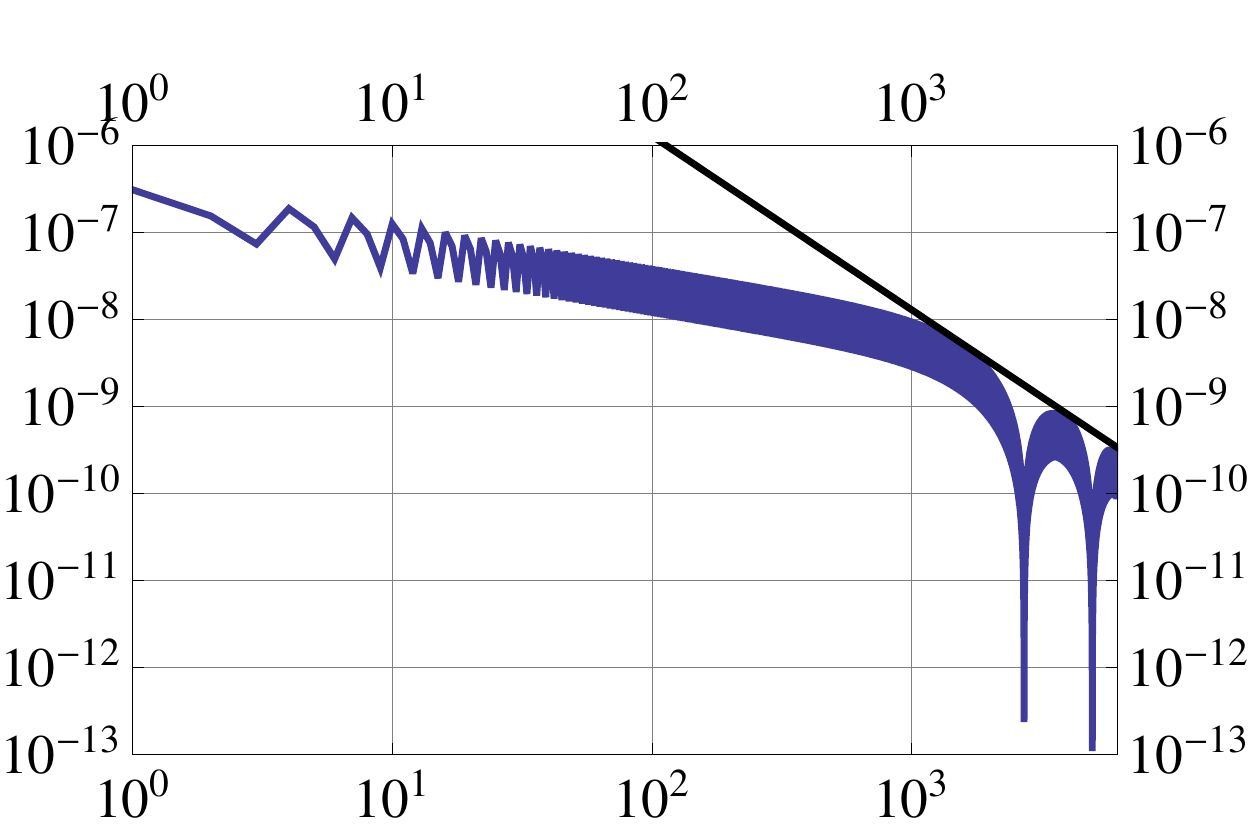}}\quad
	\subfloat[Error at $x=-99/100$, $\alpha(-99/100)=2$, $C(-99/100)=0.1245$.]
	{\includegraphics[width=0.45\textwidth]{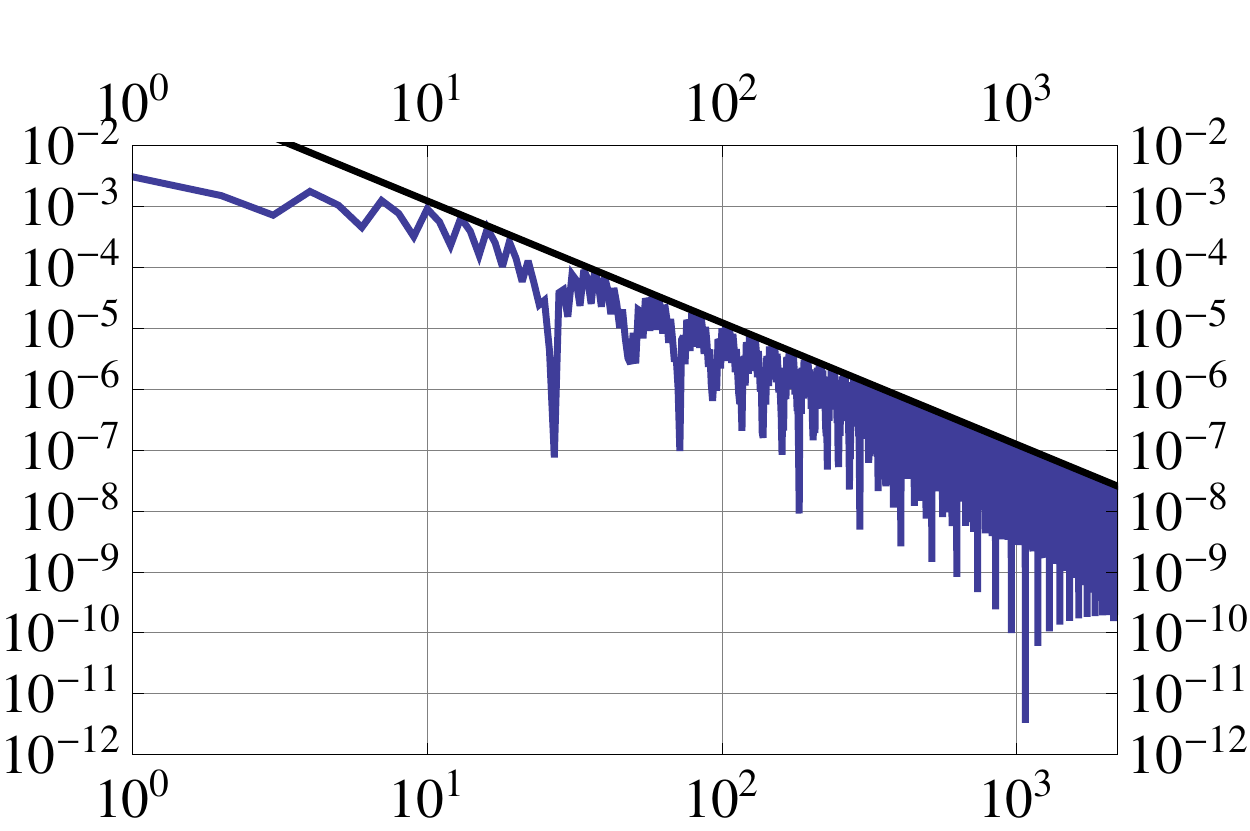}}\\
	\subfloat[Error at $x=a$, $\alpha(a)=1$, $C(a)=0.27557$.]
	{\includegraphics[width=0.45\textwidth]{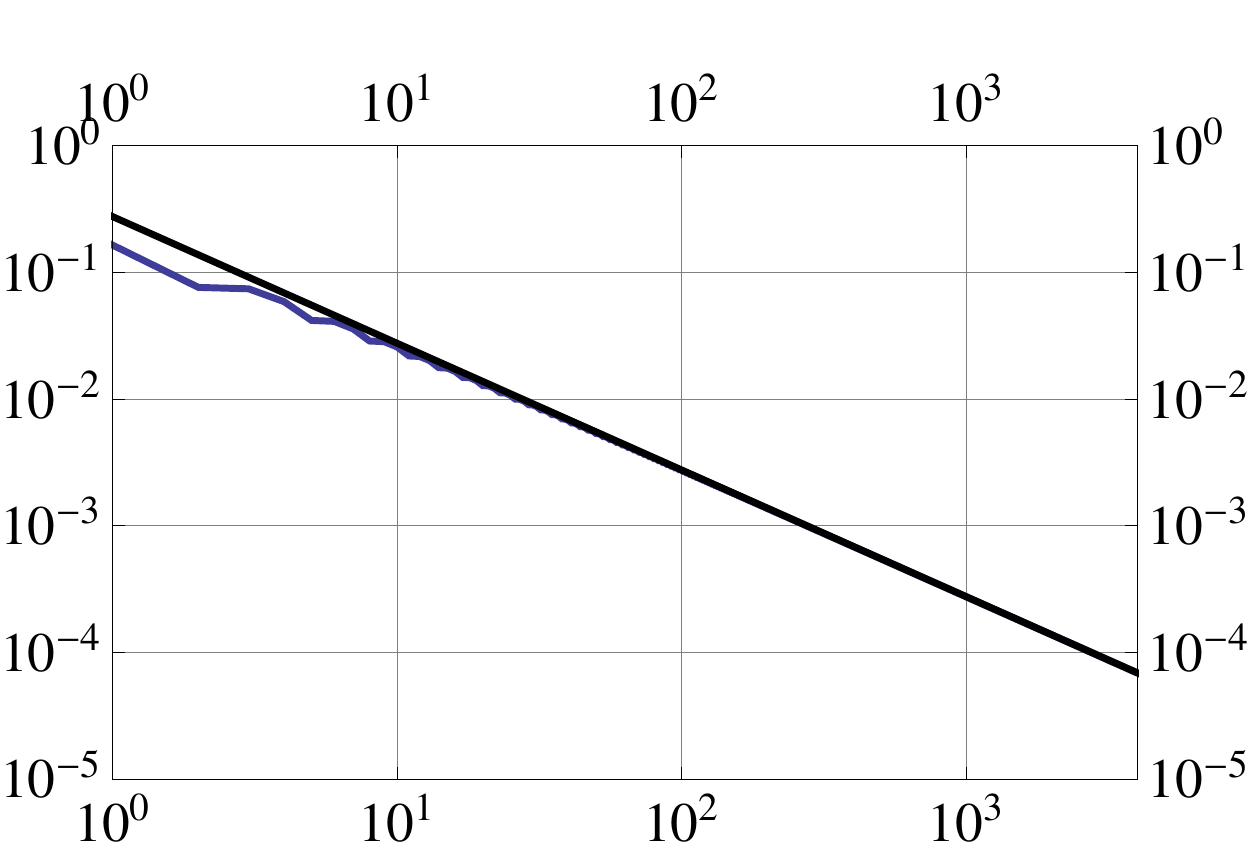}}\quad
	\subfloat[Error at $x=a+1/100$, $\alpha(a+1/100)=2$, $C(a+1/100)=0.274738$.]
	{\includegraphics[width=0.45\textwidth]{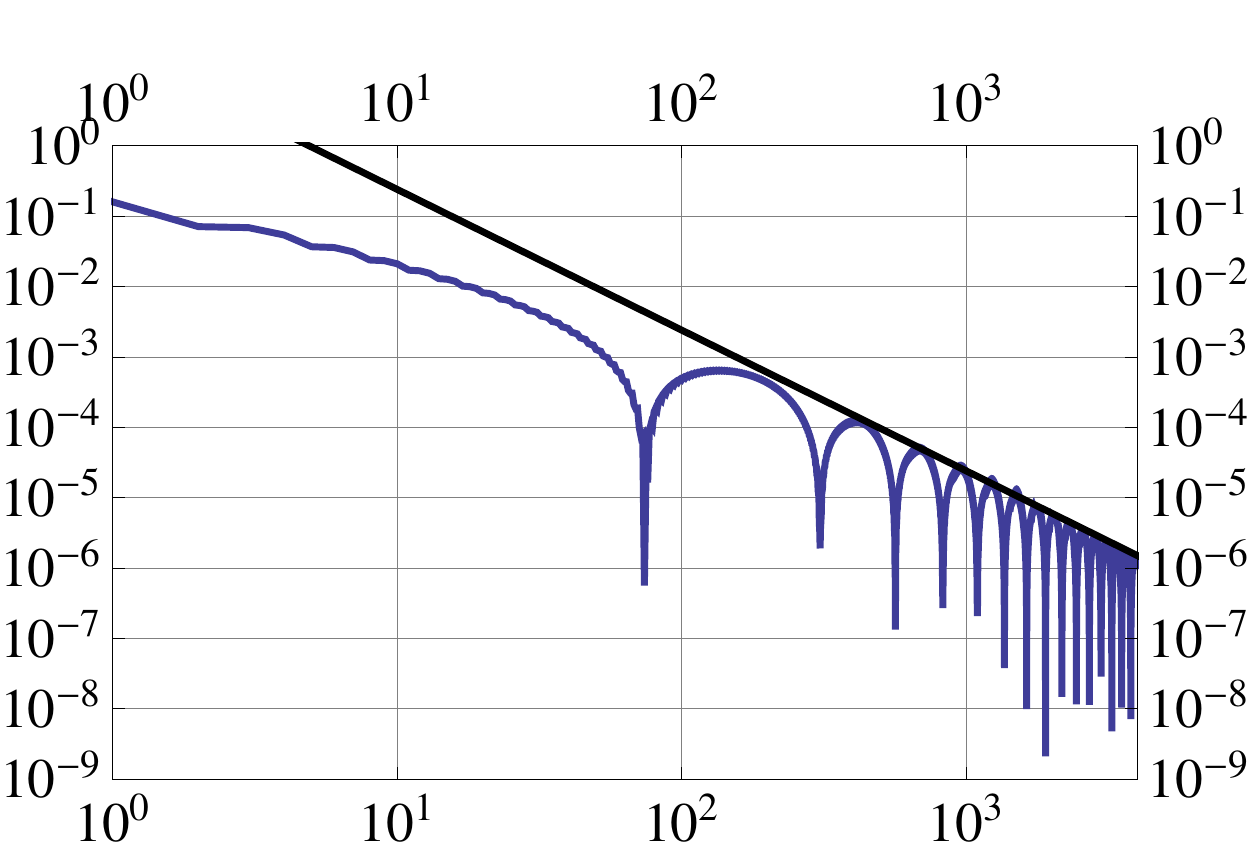}}
  \caption{$u_p(x)$: Absolute value vs polynomial order, $p=1,\ldots,2200$.}
  \label{fig:fig9}
\end{figure}
\begin{figure}[ht]
	\centering
	{\includegraphics[width=0.45\textwidth]{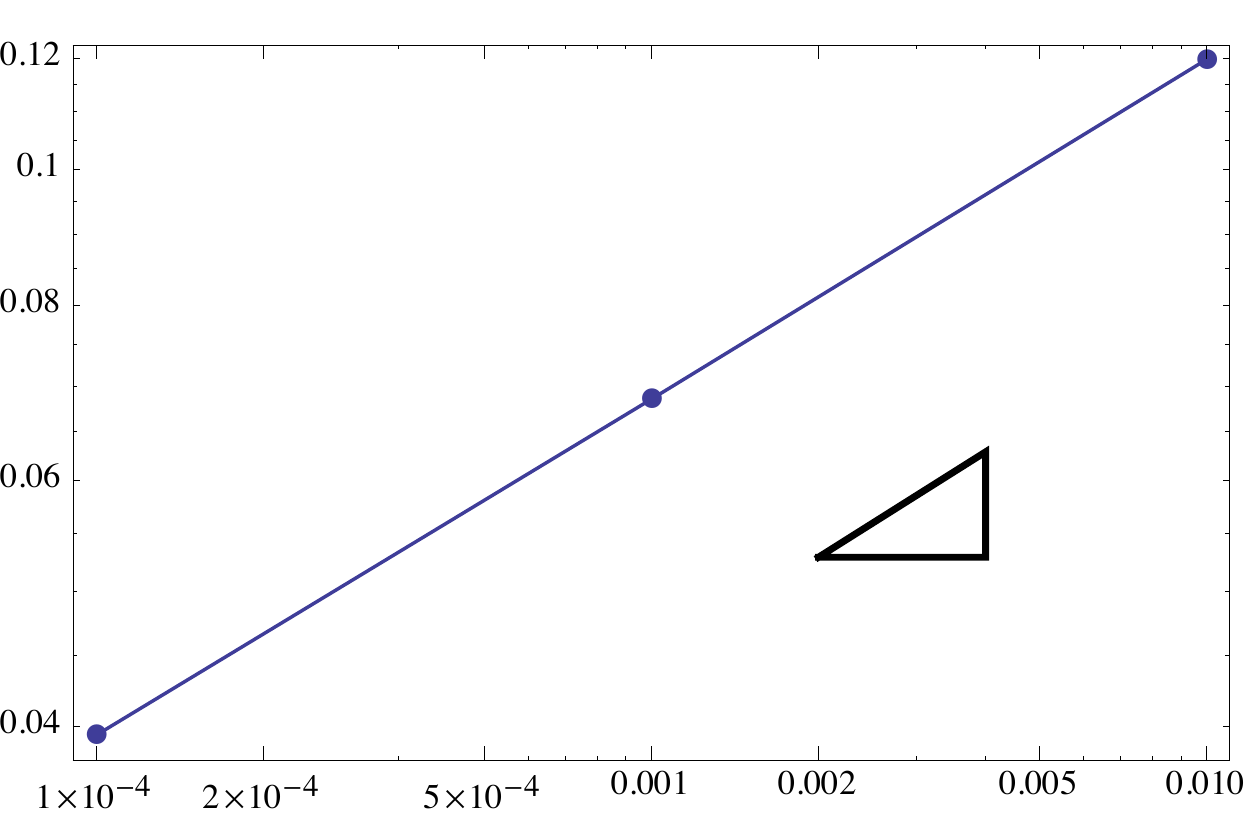}}
  \caption{$u_p(x)$: Growth of the constant $C(-1+x) \sim D\,x^{1/4}$. Value of the coefficient vs distance to the point of interest.}
  \label{fig:fig14}
\end{figure}

The results are very similar as before except that now there is no error at the
boundary. Figure~\ref{fig:fig9} shows the error for $x=-1+10^{-6},x=-1+0.01,x=a,x=a+0.01$. 
We have $\alpha (x)=2$ for $x\neq a,\alpha (a)=1$.
Comparing with the results in the
Section 6 we see the same rates.
It should be emphasized that the solution $u_{p}(x)$ is not the partial sum of the
Legendre expansion of the solution $u$.

{Let us summarize the results}
 \begin{itemize}
   \item 
     [1.] We see very similar results except in the neighborhood of the
boundary points. In the cases
when no constraint was used we had $\varepsilon ^{\prime }(-1+\xi)\sim
(\xi)^{-1/4}p^{-1}$ and $\eta (-1+\xi)\sim
(\xi)^{-1/4}p^{-2}$. In the case with the constraint we have 
$\varepsilon (-1+\xi)\sim (\xi)^{1/4}p^{-2}$, see Figure~\ref{fig:fig14}. 
\end{itemize}

\section{Generalizations}

\subsection{Legendre expansion of the function $\mid x\mid ^{\beta },-1<\beta <0$.}

\begin{figure}[ht]
	\centering
        \subfloat[Error at $x=-999/1000$, $\beta=-5/6$.]
        {\includegraphics[width=0.45\textwidth]{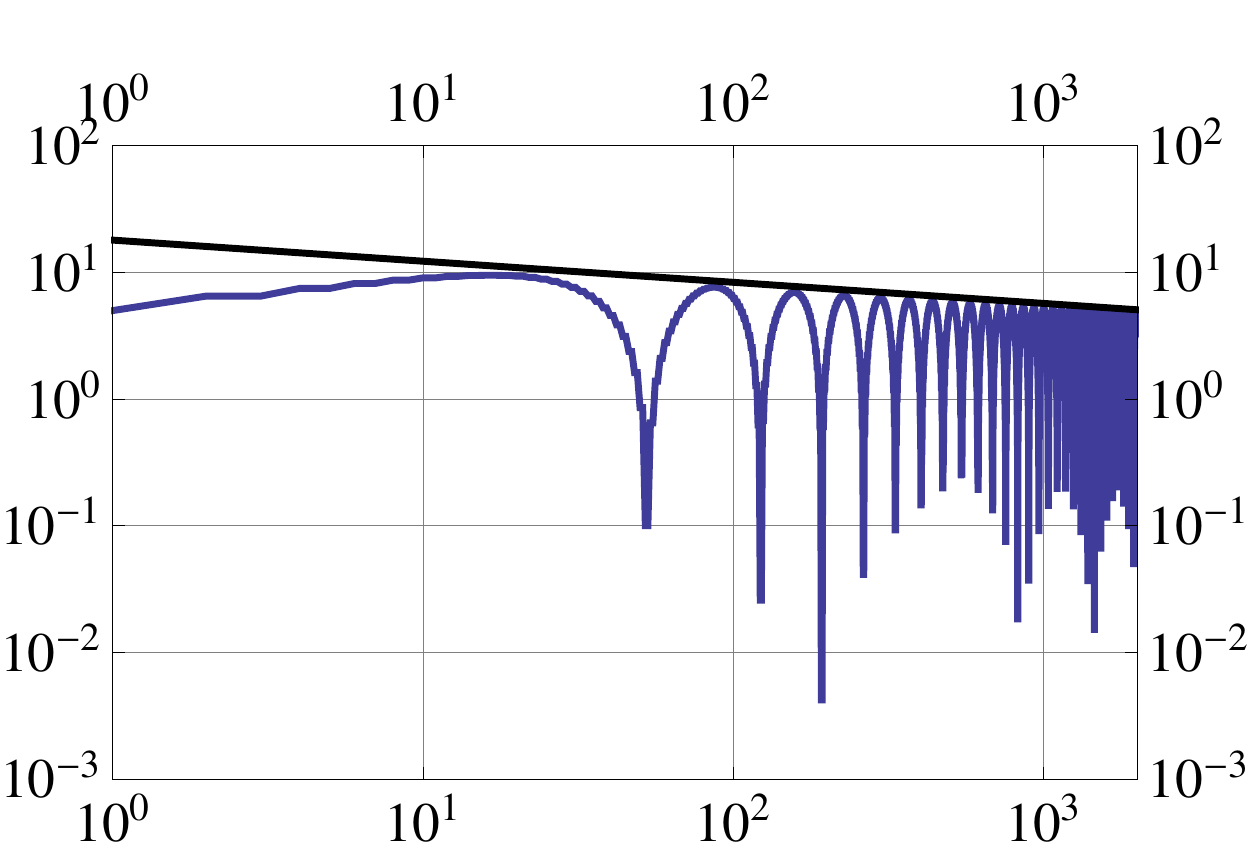}}\quad
	\subfloat[Error at $x=-99/100$, $\beta=-2/3$.]
	{\includegraphics[width=0.45\textwidth]{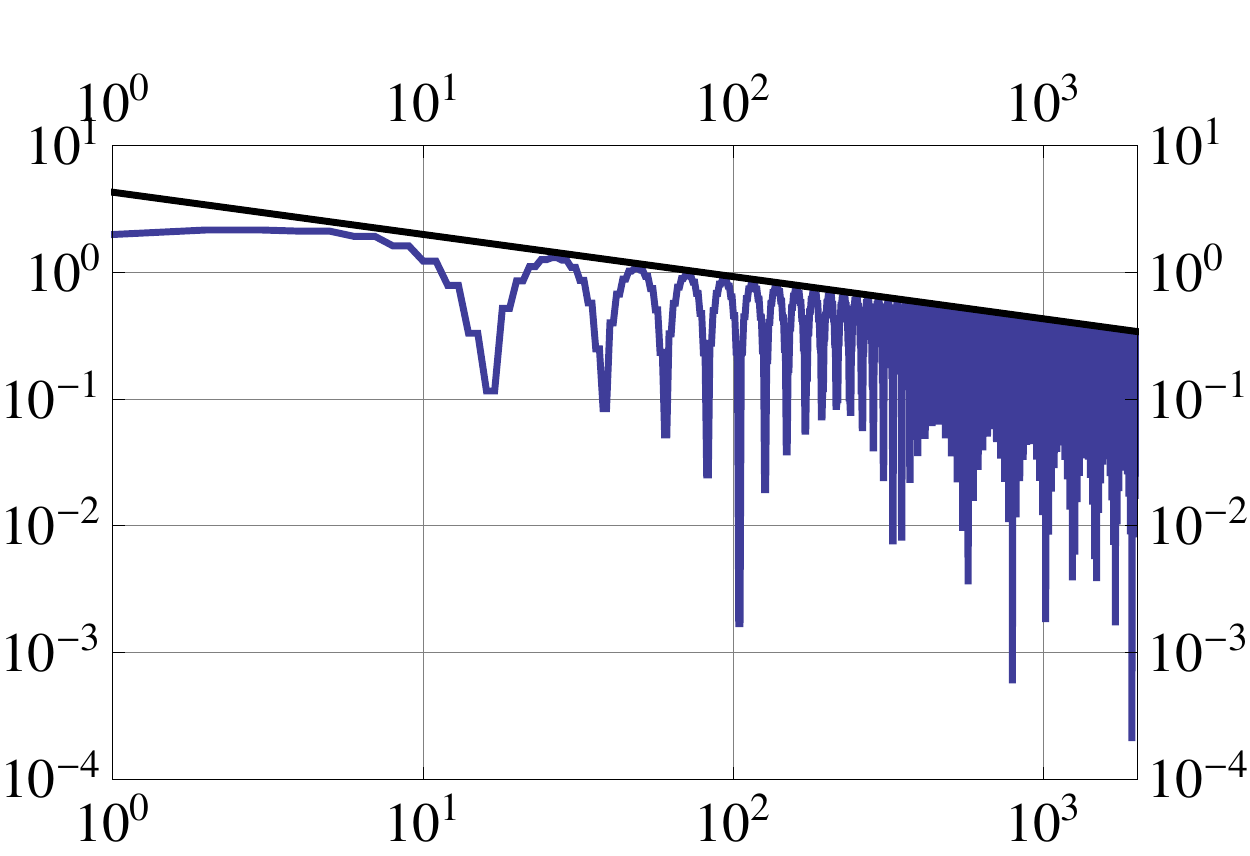}}\\
        \subfloat[Error at $x=-1/2$, $\beta=-1/2$.]
	{\includegraphics[width=0.45\textwidth]{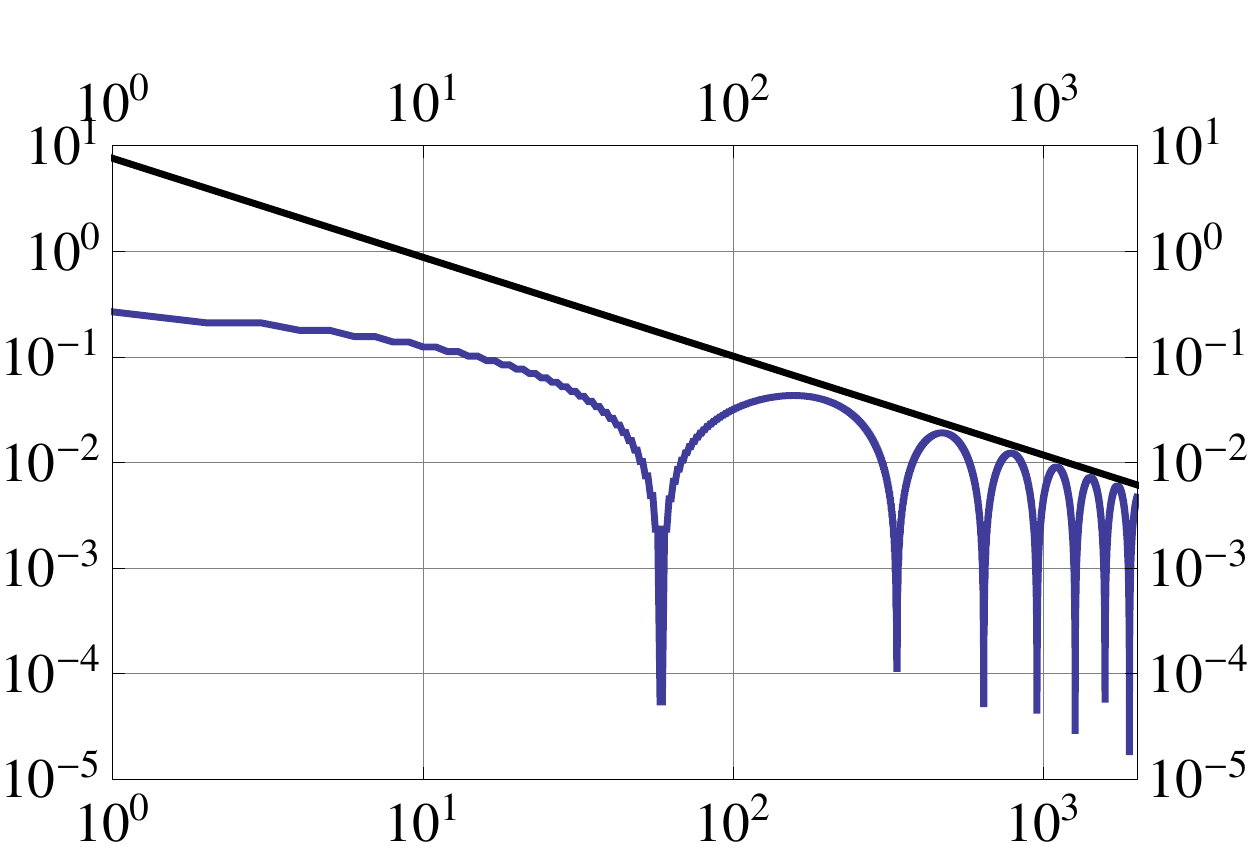}}\quad
        \subfloat[Error at $x=-1/100$, $\beta=-1/16$.]
	{\includegraphics[width=0.45\textwidth]{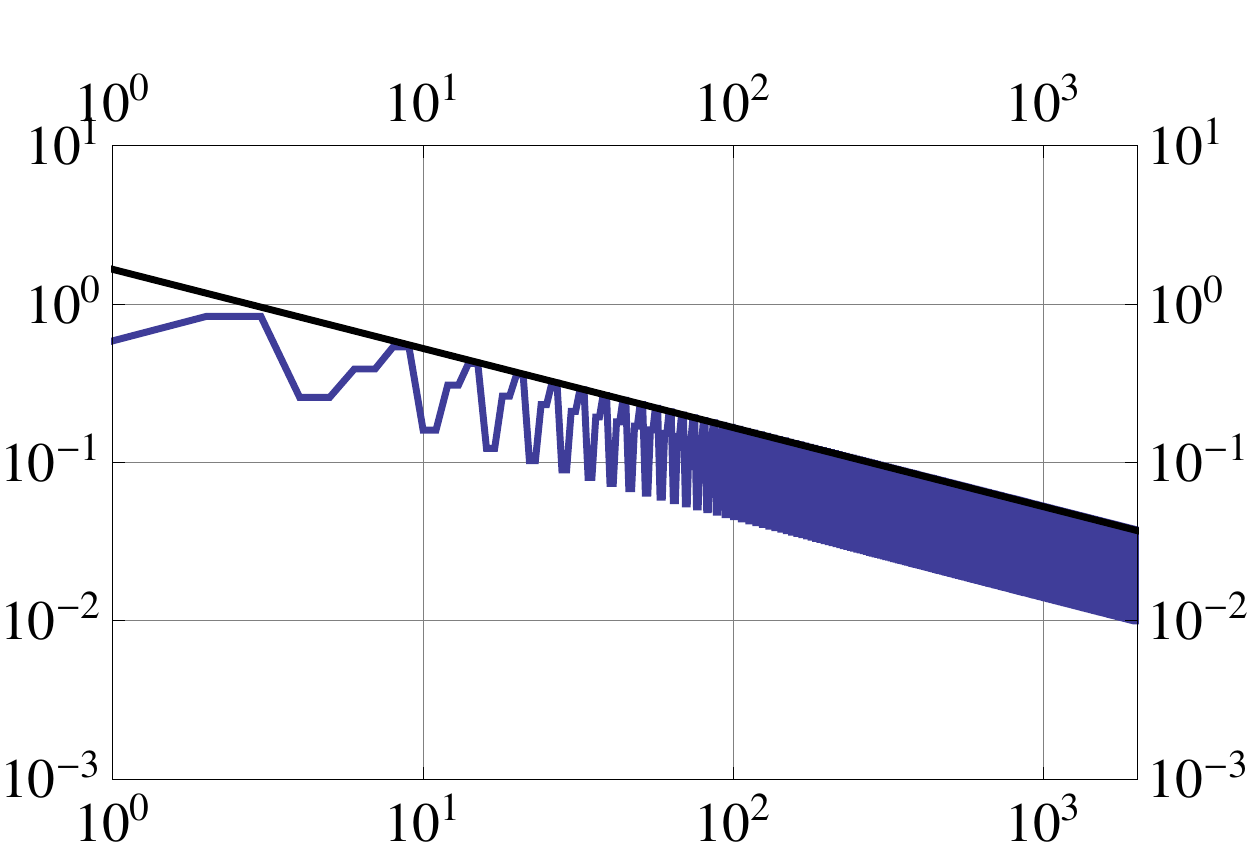}}
        \caption{Legendre expansion of the function $\mid x\mid ^{\beta}$.
      The observed rate $\alpha = 1+\beta$ in all cases. Absolute value of the error vs polynomial order, $p=1,\ldots,2200$.}
  \label{fig:alex}
\end{figure}
In the Sections 4 and 5 we addressed Legendre expansion of the
function $\mid x-a\mid$ with $a=1/2$ and $\beta =0$
resp $\beta =1$. The behavior of
the error is completely analogous for -$1<\alpha <
1$. In the case $\mid a\mid >1$ the function is
analytic and the convergence is exponential. For the theory of Legendre
expansion for analytic functions we refer to [17]. 
In this section we will address the expansion of the function $%
f(\beta,x)=$\textit{\ }$\mid x\mid ^{\beta },\ -1<\beta <0$.
Obviously we have $f(\beta ,x)\in L^{1}$ for $-1<\beta $, and $%
f(\beta ,x)\in L^{2}$ for $\beta > -1/2$. The case $\beta
=-1/2$ is a special one because $f$ is ``almost'' in $L^{2}$
and hence its integral is in $H^{1}$. This case is
important in two dimensions because the Green's function is a function of this
type. As we have already stated, the analysis of the two dimensional case is in preparation.

As shown in Figure~\ref{fig:alex}, the observed rate $\alpha = 1 + \beta$
for all $-1 < \beta < 0$.
Because of the strong singularity in $x=0$ the rate $\alpha$ has to be understood so that
\[
\mid \sum_{k=0}^{k=p}c_{k}P_{k}(0)\mid =Q(p)<Cp^{\alpha}, 
\]
with $\alpha $ negative and $C$ independent of $p$.

Then we have
\begin{itemize}
  \item
[1.] $\alpha (\pm 1)=\beta +1/2$. Hence for $\beta >-1/2$ 
we have convergence; for $\beta <-1/2$ we have divergence
with the growth $p^{\mid \alpha \mid }$. The case $\beta =-1/2$
is a special one. We see divergence with bounded partial sums.
More specifically we have $c_{2p}(\beta )P_{2p}(\pm
1)=(-1)^{p}C(\beta )p^{\beta -1/2}[1+\mathcal{O}(1/p)]$.
In general the rate $\alpha =0$ could mean one of two possibilities:
Either convergence to a wrong limit or divergence with the partial sums 
bounded.
  \item
    [2.] $\alpha (0)=\beta$. For all $-1<\beta <0$ there is
divergence bounded by $p^{\mid \beta \mid }$.
  \item
    [3.] $\alpha (x)=\beta +1$. For all $\ 0<\mid x\mid <1$ we see
convergence with the error $p^{-\mid \beta +1\mid }$. Further we have for $%
x=(-1+\xi )$, small $\xi >0$, $\mid \varepsilon (x)\mid \sim \xi
^{-1/4}p^{-(\beta +1)}$ and $Q(p)\sim \xi ^{-1}p^{-\mid \beta
+1\mid }$. 
\end{itemize}

\subsection{Legendre expansion of the function $\mid x+1 \mid ^{\beta },-1<\beta$.}
\begin{figure}[ht]
	\centering
        \subfloat[Error at $x=-1/10$, $\beta=-1/2$. Rate $\alpha = 1/2$.]
        {\includegraphics[width=0.45\textwidth]{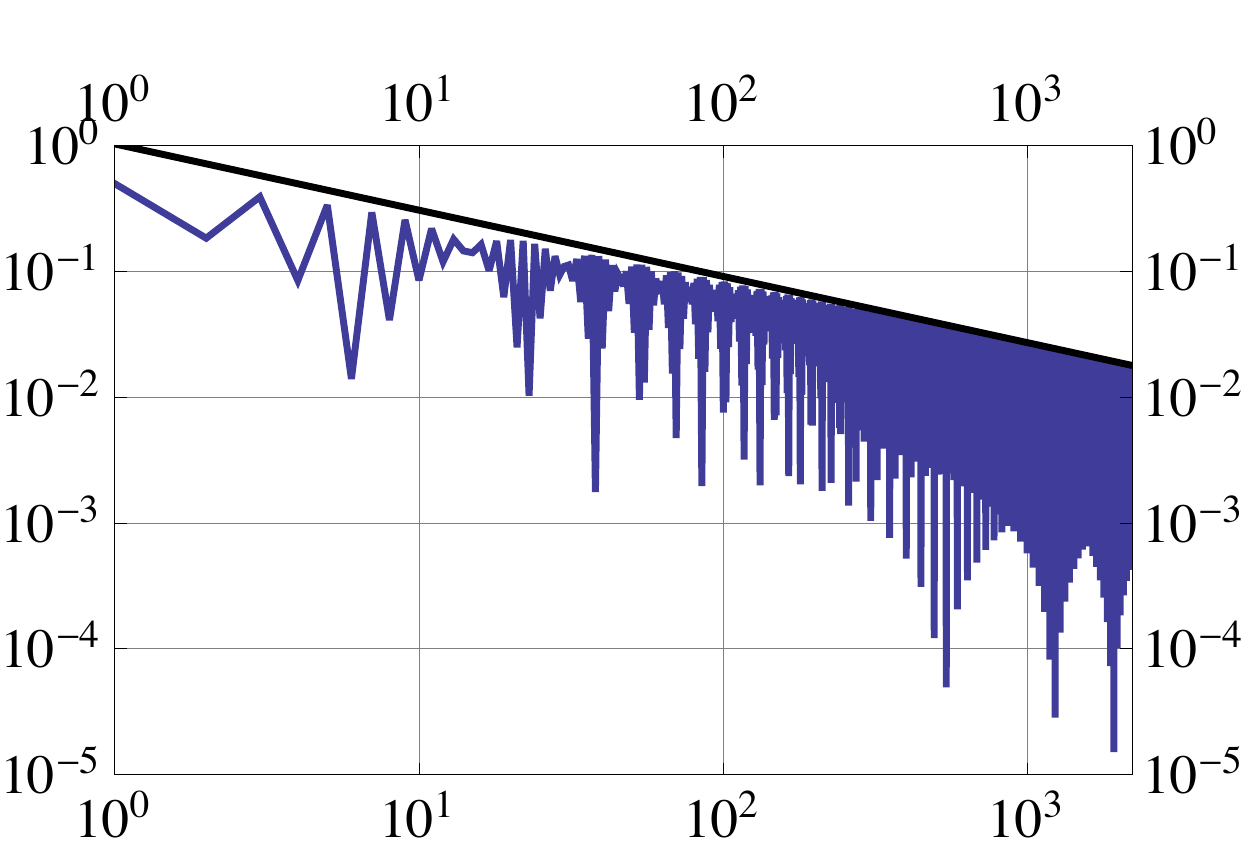}}\quad
	\subfloat[Error at $x=-1$, $\beta=1/2$. Rate $\alpha = 1$.]
	{\includegraphics[width=0.45\textwidth]{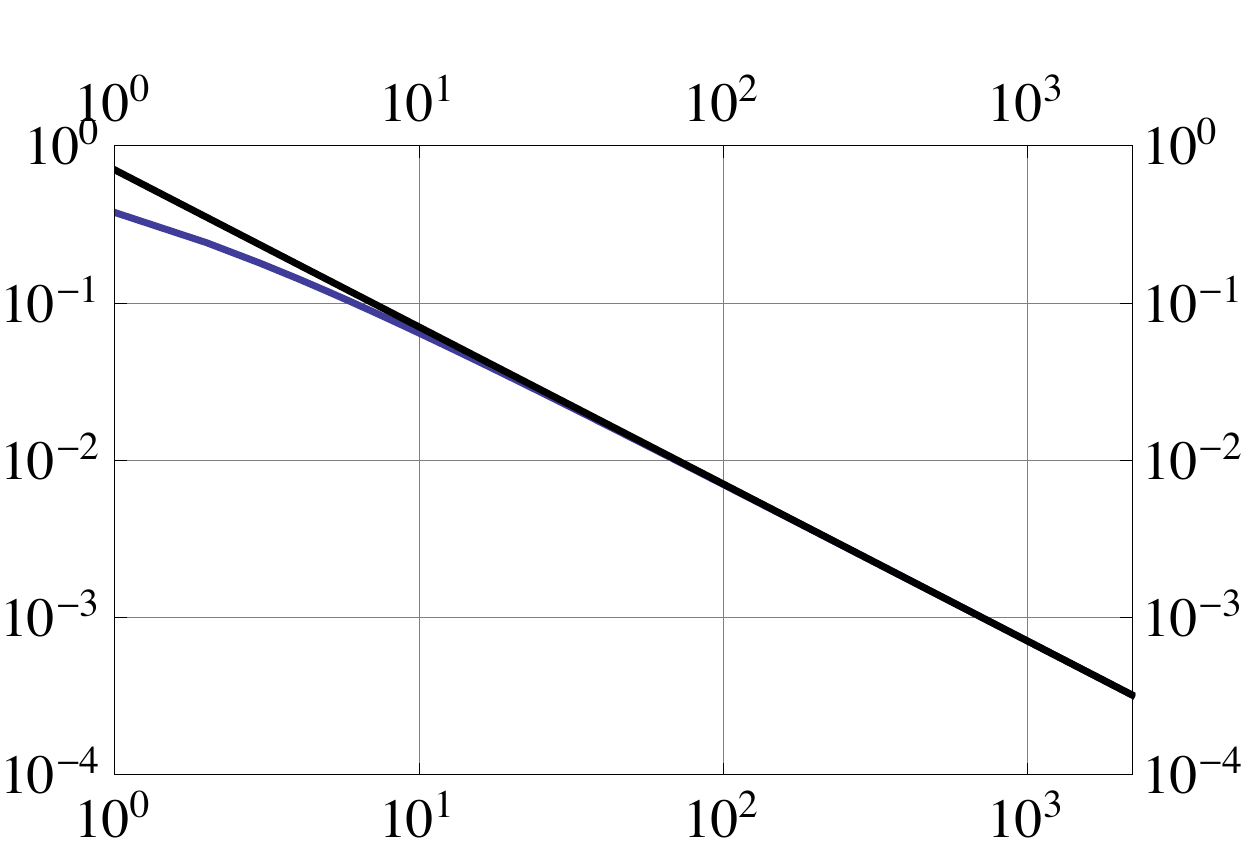}}\\
        \subfloat[Error at $x=-1/10$, $\beta=1/2$. Rate $\alpha = 5/2$.]
	{\includegraphics[width=0.45\textwidth]{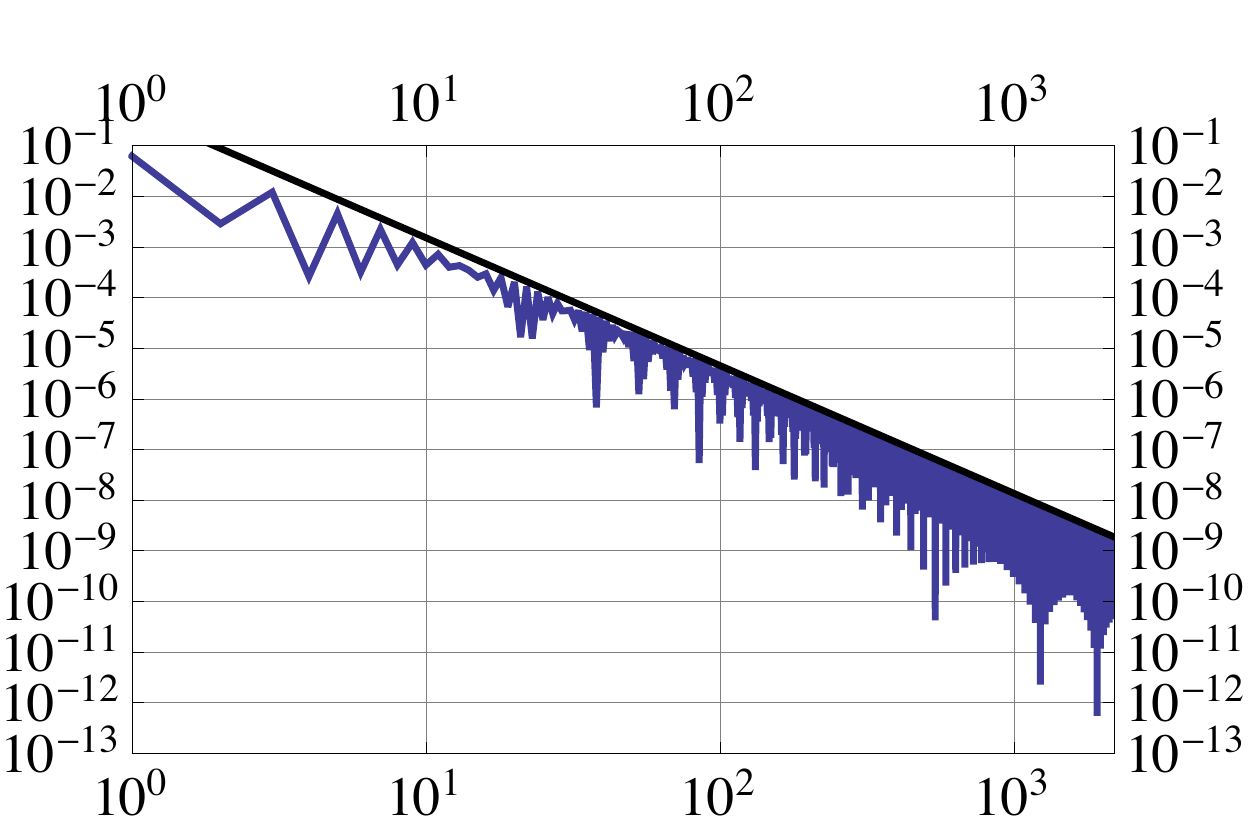}}\quad
        \subfloat[Error at $x=-1/10$, $\beta=3/2$. Rate $\alpha = 9/2$.]
	{\includegraphics[width=0.45\textwidth]{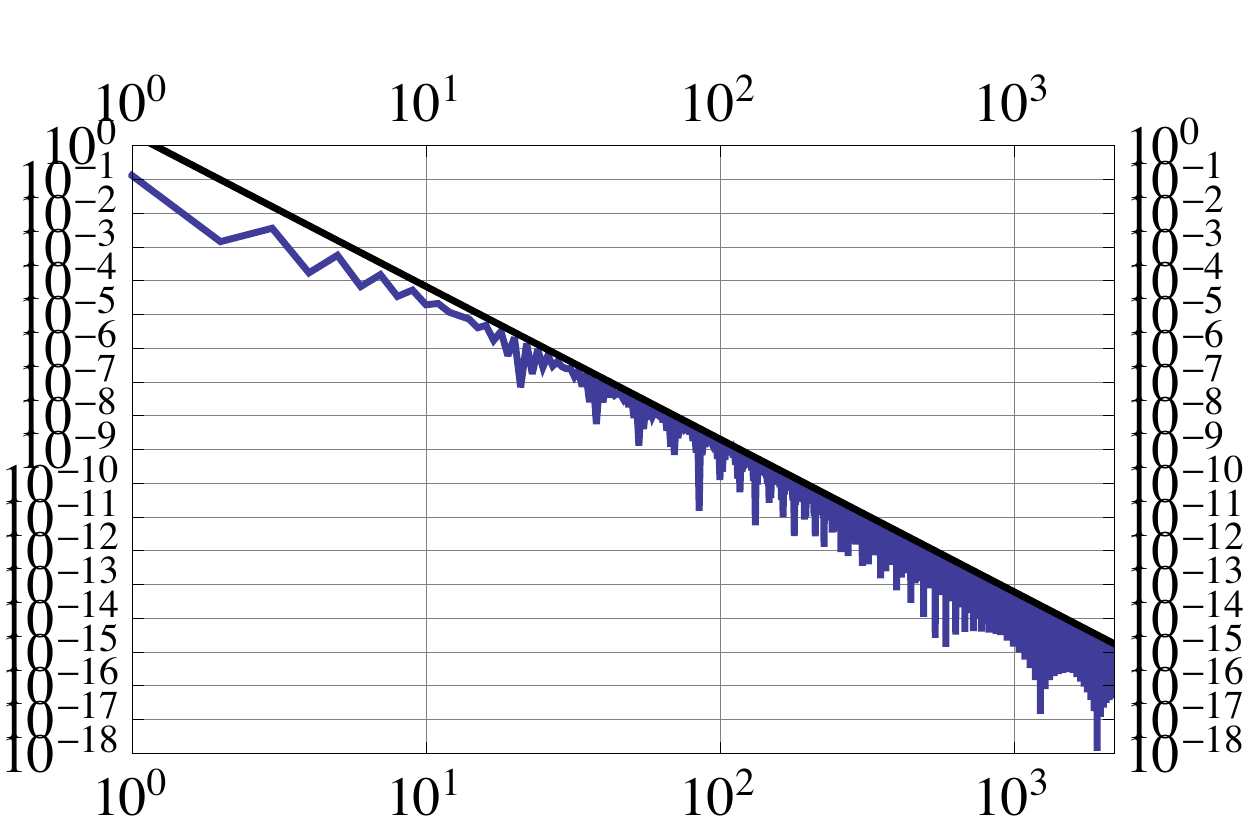}}
        \caption{Legendre expansion of the function $\mid x+1\mid ^{\beta}$. Absolute value of the error vs polynomial order, $p=1,\ldots,2200$.}
  \label{fig:alexplus}
\end{figure}

Using the conventions of the previous section we have (Figure~\ref{fig:alexplus}) 
for $\beta > -1$
\begin{itemize}
  \item [1.] $\alpha (-1)=2\beta$, naturally only for $\beta > 0$. 
  \item [2.] $\alpha (1)=2\beta+1$.
  \item [3.] $\alpha (x)=2\beta+3/2$. Further we have for 
    $x=(-1+\xi )$, small $\xi >0$, $\mid \varepsilon (x)\mid \sim \xi
^{-3/4}p^{-(2\beta+3/2)}$ and 
$x=(1-\xi )$, small $\xi >0$, $\mid \varepsilon (x)\mid \sim \xi
^{-1/4}p^{-(2\beta+3/2)}$.
\end{itemize}
Notice, that unlike in the previous examples, the coefficients of the Legendre
expansion are not known a priori.
The construction used in the computations is given in the Appendix.

\section{Summary and Conjecture}
\label{sec:hypo}
We have analyzed the behavior of the $p$-version in one
dimension on [-1,1] when the solution $u(x)$ is a piecewise analytic
function, 
\[
  u(x)=\sum_{i=1}^{m}c_{i}\mid x-a_{i}\mid ^{\beta _{i}}+v(x),\mid
a_{i}\mid <1,\beta _{i}>-1
\]
and $v(x)$ {analytic function} on
[-1,1]. We have shown that the approximate solution is the partial sum
of the Legendre expansion of $u$.

We focus on the asymptotic behavior of the $p$
partial Legendre expansion of the functions $\mid x-a\mid ^{\beta }$
leading to the algebraic convergence rate while the convergence rate for
an analytic function is exponential.

We concentrated first on the case $\beta =0,1$ in 
connection with solving a simple typical second order boundary value problem
by the $p-$version. In this context we mentioned all
theoretical estimates known to us and compared the computational results
with their theoretical predictions. Then we addressed the case for general $%
\beta $. Based on these results we formulate now the
conjecture on the error of the Legendre expansion.

\begin{conjecture}
Let $w(x)=\sum_{k=0}^{\infty }c_{k}P_{k}(x)$ be
the Legendre expansion of the function $w(x)=\mid x-a\mid ^{\beta },\mid
a\mid <1,-1<\beta $, and $w_{p}(x)=\sum_{k=0}^{p}c_{k}P_{k}(x)$ 
be its partial sum. Denote by $\varepsilon (x)=\mid w(x)-w_{p}(x)\mid $
the error. Then $\varepsilon (x)$ has the following
properties.

\begin{itemize}
  
  \item[1.] $x\in (-1,a)\cup (a,1):$ Then $\varepsilon (x)\leq
C(x)p^{-\alpha}$, $\alpha=-(\beta +1)$, where $C(x)$ independent of $p$, 
the rate $\alpha=(\beta +1)$ is optimal, i.e., it cannot be improved.

\item[2.] We have $C(-1+\xi )\leq D\xi ^{-\rho}$, $C(1-\xi )\leq
D\xi ^{-\rho}$, $\rho=1/4$, $0<\xi <\delta ,0<\delta $ with $D$ and $%
\delta $ independent of $p$ and the rate $\rho=1/4$ is optimal.

\item[3.] We have $C(a+\xi )\leq D \xi^{-\sigma}, C(a-\xi )\leq D \xi^{-\sigma}, \sigma=1, 0<\mid \xi \mid \leq
\delta ,$ with $D$ and $\delta $ independent of $p$
and the rate $\sigma=1$ is optimal.

\item[4.] $x=\pm 1;$ Then $\varepsilon (\pm 1)\leq Cp^{-\alpha}$, $\alpha=(\beta +1/2)$ 
where $C$ is independent of $p$. The rate $%
\alpha=(\beta +1/2)$ is optimal. For $\beta <-1/2$ 
there is divergence and for $\beta =-1/2$ there is convergence to
a  limit which is not $w(\pm 1)$. In general the rate $\alpha=0$ 
will be understood as indicating a bounded sequence.

\item[5.] $x=a:$ Then $\varepsilon (a)\leq Cp^{-\alpha }$, where $C$
is independent of $p$. The rate $\alpha=\beta $ is
optimal. The rate $\alpha =0$ is special because the function $w$ 
is then discontinuous. Denoting $\overline{w}(a)=\frac{1}{2}%
(w(a+0)+w(a-0))$. Then there is convergence to $\overline{w}$ with the
rate $\alpha=1$.
\end{itemize}
\end{conjecture}

Some comments:
Statement 2 is related to the boundary layer because the rates in $%
x=-1$ and $x>-1$ are different. 
Statement 3 is related to the Gibbs phenomenon which
appears for all $\beta$.
Note that in the above statements the term $\lg p$ is not present.

We have seen typical features of the $p-$version for
piecewise analytic solutions.
\begin{itemize}
  
  \item[1.] The statements in the above theorems are not covered by the
available mathematical theory. They are conjectures  based on the careful
computations and their generalization.

\item[2.] The errors are highly oscillatory
and the pattern is different for different values of $x$.

\item[3.] The preasymptotic range is large and the practical computations
likely outside the asymptotic range.

\item[4.] The error behavior in the energy norm and the $L^{2}$ 
is well covered by the theory, is not  oscillatory, and the preasymptotic
range is much shorter.

\item[5.] When the solution is very smooth, precisely analytic on the entire
domain [-1.1], the convergence is exponential.

\item[6.] The behavior of the $p$-version has some but not all characteristics
  in higher dimensions. 
  We shall address these issues in the future.
\end{itemize}

\begin{acknowledgement}
  Authors would like to thank Prof P. Nevai, Ohio State University, and Prof G. Mastroianni,
  University Basilicata, Potenza, Italy, for mentioning us some known theoretical results.
\end{acknowledgement}

\appendix
\section%
{Legendre Coefficients of $\mid x+1 \mid ^{\beta },-1<\beta$.}
Consider the identity
\begin{equation}\label{eq:int}
  \begin{split}
  I_k=\int_{-1}^{1}  \mid x+1 \mid ^{\beta }x^k dx = \frac{\, _2F_1(k+1,-\beta ;k+2;-1)}{k+1}+ \\
  +\frac{(-1)^k \Gamma (\beta +1) \Gamma (k+1)}{\Gamma (k+\beta +2)},\quad \beta > -1,\ k \in \mathbb{Z},k\geq 0,
  \end{split}
\end{equation}
where $\, _2F_1(a,b;c,z)$ denotes the hypergeometric function.
Notice that every term converges for $-1<\beta$, since $\, _2F_1(a,b;c,z)$
converges for $z=\pm 1$, if $c > a+b$, that is, $k+2 > k+1-\beta$.

Expanding the Legendre polynomial $P_k(x)=\sum_{m=0}^k a_m x^m$, the standard expansion
\[
  c_k = \int_{-1}^{1} \mid x+1 \mid ^{\beta }P_k(x) dx\, (2k+1)/2,
\]
becomes using the identity of (\ref{eq:int})
\begin{equation}
 c_k = (\sum_{m=0}^k a_m I_m) (2k+1)/2.
\end{equation}

\bibliographystyle{model1-num-names}
\bibliography{Paper1}

\begin{thebibliography}{26}
\expandafter\ifx\csname natexlab\endcsname\relax\def\natexlab#1{#1}\fi
\providecommand{\bibinfo}[2]{#2}
\ifx\xfnm\relax \def\xfnm[#1]{\unskip,\space#1}\fi
\bibitem[{Babu\v{s}ka et~al.(1981)Babu\v{s}ka, Szabo, and Katz}]{c1}
\bibinfo{author}{I.~Babu\v{s}ka}, \bibinfo{author}{B.~Szabo},
  \bibinfo{author}{I.~Katz},
\newblock \bibinfo{title}{The $p$-version of the finite element method},
\newblock \bibinfo{journal}{SIAM J Numer Anal} \bibinfo{volume}{18}
  (\bibinfo{year}{1981}) \bibinfo{pages}{515--545}.
\bibitem[{Babu\v{s}ka and Dorr(1981)}]{c2}
\bibinfo{author}{I.~Babu\v{s}ka}, \bibinfo{author}{M.~R. Dorr},
\newblock \bibinfo{title}{Error estimate for the combined $h$- and $p$-
  versions of the finite element method,},
\newblock \bibinfo{journal}{Numer. Math.} \bibinfo{volume}{37}
  (\bibinfo{year}{1981}) \bibinfo{pages}{257--277}.
\bibitem[{Oden(1991)}]{c3}
\bibinfo{author}{J.~T. Oden},
\newblock \bibinfo{title}{Finite elements: An introduction},
\newblock in: \bibinfo{editor}{P.~G. Ciarlet}, \bibinfo{editor}{J.~L. Lions}
  (Eds.), \bibinfo{booktitle}{Handbook of Numerical Analysis}, volume
  \bibinfo{volume}{II, Finite Element Methods, Part I},
  \bibinfo{publisher}{North Holland}, \bibinfo{year}{1991}, pp.
  \bibinfo{pages}{3--17}.
\bibitem[{Babu\v{s}ka(1994)}]{c4}
\bibinfo{author}{I.~Babu\v{s}ka},
\newblock \bibinfo{title}{Courant element: Before and after in finite element
  methods},
\newblock in: \bibinfo{editor}{M.~K. Krizek},
  \bibinfo{editor}{P.~Neittaanmaki}, \bibinfo{editor}{R.~Stenberg} (Eds.),
  \bibinfo{booktitle}{Fifty years of Courant method}, volume
  \bibinfo{volume}{164} of \textit{\bibinfo{series}{Lecture Notes in Pure and
  Applied Mathematics}}, \bibinfo{publisher}{Marcel Dekker},
  \bibinfo{year}{1994}, pp. \bibinfo{pages}{37--51}.
\bibitem[{Babu\v{s}ka and Aziz(1973)}]{c5}
\bibinfo{author}{I.~Babu\v{s}ka}, \bibinfo{author}{A.~Aziz},
\newblock \bibinfo{title}{Survey lectures on the mathematical foundations of
  finite element method},
\newblock in: \bibinfo{editor}{A.~Aziz} (Ed.), \bibinfo{booktitle}{The
  Mathematical Foundations of the Finite Element Methods with Applications to
  Partial Differential Equations,}, \bibinfo{publisher}{Academic Press},
  \bibinfo{year}{1973}.
\bibitem[{Ciarlet(1978)}]{c6}
\bibinfo{author}{P.~Ciarlet}, \bibinfo{title}{The Finite Element Method for
  Elliptic Problems}, \bibinfo{publisher}{North Holland}, \bibinfo{year}{1978}.
\bibitem[{Wahlbin(1991)}]{c7}
\bibinfo{author}{L.~Wahlbin},
\newblock \bibinfo{title}{Local behavior in finite element method},
\newblock in: \bibinfo{editor}{P.G.Ciarlet}, \bibinfo{editor}{J.~L. Lions}
  (Eds.), \bibinfo{booktitle}{Handbook of Numerical Analysis}, volume
  \bibinfo{volume}{I, Finite Element Methods, Part I},
  \bibinfo{publisher}{North Holland}, \bibinfo{year}{1991}, pp.
  \bibinfo{pages}{353--523}.
\bibitem[{Szabo and Babu\v{s}ka(1991)}]{bc1}
\bibinfo{author}{B.~Szabo}, \bibinfo{author}{I.~Babu\v{s}ka},
  \bibinfo{title}{Finite {E}lement {A}nalysis}, \bibinfo{publisher}{Wiley},
  \bibinfo{year}{1991}.
\bibitem[{Szabo and Babu\v{s}ka(2011)}]{bc2}
\bibinfo{author}{B.~Szabo}, \bibinfo{author}{I.~Babu\v{s}ka},
  \bibinfo{title}{Introduction to {F}inite {E}lement {A}nalysis},
  \bibinfo{publisher}{Wiley}, \bibinfo{year}{2011}.
\bibitem[{Demkowicz(2006)}]{c20}
\bibinfo{author}{L.~Demkowicz}, \bibinfo{title}{Computing with $hp$-ADAPTIVE
  FINITE ELEMENTS: Volume 1 One and Two Dimensional Elliptic and Maxwell
  Problems}, \bibinfo{publisher}{Chapman and Hall}, \bibinfo{year}{2006}.
\bibitem[{Demkowicz(2007)}]{c21}
\bibinfo{author}{L.~Demkowicz}, \bibinfo{title}{Computing with hp-ADAPTIVE
  FINITE ELEMENTS: Volume II Frontiers: Three Dimensional Elliptic and Maxwell
  Problems with Applications}, \bibinfo{publisher}{Chapman and Hall},
  \bibinfo{year}{2007}.
\bibitem[{Suetin(1978)}]{c9}
\bibinfo{author}{P.~K. Suetin}, \bibinfo{title}{Classical Orthogonal
  Polynomials (in Russian)}, \bibinfo{publisher}{Nauka}, \bibinfo{year}{1978}.
\bibitem[{Jackson(1930)}]{c10}
\bibinfo{author}{D.~Jackson}, \bibinfo{title}{The Theory of Approximation},
  \bibinfo{publisher}{American Mathematical Society}, \bibinfo{year}{1930}.
\bibitem[{Szego(1939)}]{c11}
\bibinfo{author}{C.~Szego}, \bibinfo{title}{Orthogonal Polynomials},
  \bibinfo{publisher}{American Mathematical Society}, \bibinfo{year}{1939}.
\bibitem[{Borwein and Erdelyi(1975)}]{c12}
\bibinfo{author}{P.~Borwein}, \bibinfo{author}{V.~Erdelyi},
  \bibinfo{title}{Polynomials and Polynomial Inequalities},
  \bibinfo{publisher}{Springer}, \bibinfo{year}{1975}.
\bibitem[{D.Gottlieb and Shu(1997)}]{c13}
\bibinfo{author}{D.Gottlieb}, \bibinfo{author}{C.-W. Shu},
\newblock \bibinfo{title}{Gibbs phenomenon and its resolution},
\newblock \bibinfo{journal}{SIAM Review} \bibinfo{volume}{39}
  (\bibinfo{year}{1997}) \bibinfo{pages}{644--668}.
\bibitem[{Hewitt and Hewitt(1979)}]{c14}
\bibinfo{author}{E.~Hewitt}, \bibinfo{author}{R.~E. Hewitt},
\newblock \bibinfo{title}{The {G}ibbs-{W}ilbraham phenomenon. an episode in
  {F}ourier analysis},
\newblock \bibinfo{journal}{Archive for History of Exact Sciences}
  \bibinfo{volume}{21} (\bibinfo{year}{1979}) \bibinfo{pages}{129--160}.
\bibitem[{Benardi and Maday(1997)}]{c8}
\bibinfo{author}{C.~Benardi}, \bibinfo{author}{Y.~Maday},
\newblock \bibinfo{title}{Spectral methods},
\newblock in: \bibinfo{editor}{P.~G. Ciarlet}, \bibinfo{editor}{J.~L. Lions}
  (Eds.), \bibinfo{booktitle}{Handbook of Numerical Analysis}, volume
  \bibinfo{volume}{V, Techniques of Scientific Computing},
  \bibinfo{publisher}{Elsevier}, \bibinfo{year}{1997}, pp.
  \bibinfo{pages}{209--487}.
\bibitem[{Shen et~al.(2011)Shen, Tang, and Wang}]{c15}
\bibinfo{author}{J.~Shen}, \bibinfo{author}{T.~Tang}, \bibinfo{author}{L.-L.
  Wang}, \bibinfo{title}{Spectral Methods, Algorithms, Analysis and
  Applications}, \bibinfo{publisher}{Springer}, \bibinfo{year}{2011}.
\bibitem[{Mastroianni and Milovanovi\'{c}(2008)}]{cm2}
\bibinfo{author}{G.~Mastroianni}, \bibinfo{author}{G.~V. Milovanovi\'{c}},
  \bibinfo{title}{Interpolation processes. Basic theory and applications},
  Springer Monographs in Mathematics, \bibinfo{publisher}{Springer--Verlag},
  \bibinfo{address}{Berlin}, \bibinfo{year}{2008}.
\bibitem[{Sundar et~al.(2013)Sundar, Stadlerand, and Biros}]{bc3}
\bibinfo{author}{H.~Sundar}, \bibinfo{author}{G.~Stadlerand},
  \bibinfo{author}{G.~Biros}, \bibinfo{title}{Comparison of multigrid
  algorithms for high-order continuous finite element discretizations},
  \bibinfo{year}{2013}.
\bibitem[{Bojanic and Vuilleumier(1981)}]{c16}
\bibinfo{author}{R.~Bojanic}, \bibinfo{author}{M.~Vuilleumier},
\newblock \bibinfo{title}{On the rate of convergence of {F}ourier-{L}egendre
  series of functions of bounded variation},
\newblock \bibinfo{journal}{Journal of Approximation Theory}
  \bibinfo{volume}{31} (\bibinfo{year}{1981}) \bibinfo{pages}{67--79}.
\bibitem[{Mastroianni(2013)}]{cM}
\bibinfo{author}{G.~Mastroianni}, \bibinfo{title}{Private communication},
  \bibinfo{year}{2013}.
\bibitem[{Kudromonov(2006)}]{c17}
\bibinfo{author}{D.~Kudromonov},
\newblock \bibinfo{title}{Convergence estimate of the {F}ourier-{L}egendre
  series of functions with bounded variation (in russian)},
\newblock \bibinfo{journal}{Izvestia vyshich ucebnych zavedenij, Matematika}
  \bibinfo{volume}{7(530)} (\bibinfo{year}{2006}) \bibinfo{pages}{34--45}.
\bibitem[{Wang and Xiang(2012)}]{c18}
\bibinfo{author}{H.~Wang}, \bibinfo{author}{S.~Xiang},
\newblock \bibinfo{title}{On the convergence rates of {L}egendre
  approximation},
\newblock \bibinfo{journal}{Mathematics of Computation} \bibinfo{volume}{81}
  (\bibinfo{year}{2012}) \bibinfo{pages}{861--877}.
\bibitem[{Lebedev(1965)}]{c9a}
\bibinfo{author}{N.~N. Lebedev}, \bibinfo{title}{Special functions and their
  applications}, \bibinfo{publisher}{Prentice Hall}, \bibinfo{year}{1965}.

\end{thebibliography}
\end{document}